\documentclass{article}

\usepackage{amssymb}
\usepackage{amsmath}
\usepackage{amscd}
\usepackage{amsthm}
\usepackage[dvips]{color}
\usepackage{graphicx}
\usepackage{pb-diagram}
\usepackage{color}
\usepackage{epsfig}

\newcommand{\C} {\mathbb C}

\newcommand{\h} {\hat}
\newcommand{\mc}{\mathcal}

\newcommand{\cantor}{\{0,1\}^{\mathbb{N}}}

\newcommand{\inver}{\underset{\longleftarrow}{\lim}}
\newcommand{\fiber}{\pi^{-1}}

\newcommand{\ra}{\rightarrow}

\newcommand{\cone}{Con\,(\mc{S}^1)}
\newcommand{\sol}{\mathcal{S}^1}
\newcommand{\NN}{\mathcal{N}}
\newcommand{\RR}{\mathcal{R}}

\newtheorem{theorem}{Theorem}
\newtheorem{lemma}[theorem]{Lemma}
\newtheorem{corollary}[theorem]{Corollary}
\newtheorem{proposition}[theorem]{Proposition}
\newtheorem*{definition}{Definition}
\newtheorem*{thm}{Main Theorem}

\title{On the classification of laminations associated to quadratic polynomials}
\author{Carlos Cabrera}
\begin{document}
\maketitle
\begin{abstract}
Given any rational map $f$, there is a lamination by Riemann
surfaces associated to $f$. Such laminations were constructed in
general by Lyubich and Minsky. In this paper, we classify
laminations associated to quadratic polynomials with periodic
critical point. In particular, we prove that the topology of such
laminations determines the combinatorics of the parameter. We also
describe the topology of laminations associated to other types of
quadratic polynomials.

\end{abstract}

\tableofcontents

\section{Introduction}

The ``natural extension'' $\NN_f$ (or the ``inverse limit'') of a
rational map $f$  is a very interesting object whose topology and
geometry reflects in an intricate way the dynamics of $f$ (see
Lyubich and Minsky in \cite{LM}). In this paper,  we study the
relation between the topology of $\NN_f$ and the dynamics of $f$,
focusing on the case of quadratic polynomials $f_c: z\mapsto
z^2+c$.

The natural extension $\NN_f$  contains the ``regular leaf space''
$\RR_f$ whose connected components (``leaves'') are endowed with a
natural conformal structure. When the dynamics of $f$ is simple,
the corresponding regular leaf space has a lamination structure;
that is, there is an atlas of charts, such that the image of every
chart is the product of a disk times a Cantor set. Moreover, the
leaves of $\RR_f$ are simply connected Riemann surfaces
conformally isomorphic to the complex plane. This is the case for
hyperbolic rational maps, where $\RR_f$ is obtained from $\NN_f$
by removing finitely many points. In general, associated to any
rational map, Lyubich and Minsky's construction provides us with
an \textit{orbifold affine lamination} which, in turn, admits an
extension to a \textit{3-dimensional hyperbolic lamination}.

The main goal of this paper is to prove that, for hyperbolic
quadratic polynomials, the topology of the lamination determines
the combinatorics of the corresponding parameter. More precisely:

\begin{thm}
If $h:\mc{R}_c\ra \mc{R}_{c'}$ is an orientation preserving
homeomorphism between the regular leaf spaces associated to the
superattracting parameters $c$ and $c'$, then $c=c'$. \
\end{thm}
\medskip
Let us now give a more detailed description of the contents of the
paper.

In Section 2, we summarize the necessary background in basic
holomorphic dynamics. We assume the reader is familiar with the
subject, and we highlight only the facts that are important for
understanding the structure of the associated laminations. Most of
this theory is readily available in various surveys in complex
dynamics, e.g.,\cite{Be},\cite{CG}, \cite{DH}, \cite{L} and
\cite{Mdyn}.

From the lamination point of view, the classical theory of
linearization around periodic points provides us with suitable
uniformizations of the corresponding periodic leaves. We will make
use of this application of K\"onigs and B\"ottcher coordinates
(near repelling and superattracting points). This approach has
proven useful in the parabolic case (see Tomoki Kawahira
\cite{Tom}).

Although we are focused on the superattracting case, most of our
results are valid over a broader class of parameters, specifically
those for which the critical point is non-recurrent and has no
parabolics. According to \cite{LM}, a rational map has
non-recurrent critical points and no parabolics if and only if the
corresponding 3-lamination is convex cocompact. Because of this
property, we call this class of quadratic polynomials, or the
corresponding parameters, \textit{convex cocompact}. We describe
some basic dynamical features of convex cocompact polynomials (see
\cite{CJY}
 and \cite{Man}).

There are several combinatorial models that describe
superattracting quadratic polynomials. For a complete discussion
of them see the paper of Henk Bruin and Dierk Schleicher
\cite{BS}. We describe some of these models, for which we have
found topological analogues in the regular leaf space. One of the
most informative is the model of \textit{ray portraits} as
presented by John Milnor in \cite{Mext}. For the proof of our main
theorem, we will see that the topology of the  affine lamination
associated with a superattracting quadratic polynomial determines
the ray portrait of the corresponding parameter. The interested
reader can also find more about the combinatorics of
postcritically finite quadratic polynomials in \cite{BS},
\cite{D},\cite{DH},\cite{HFB}, \cite{LS}, \cite{Mext}, and
\cite{Poi}.

\medskip
In Section 3, we discuss basic definitions and properties of
inverse limits. A classical example of an inverse limit is the
{\it dyadic solenoid} $\sol$, associated to the map $f_0: z\mapsto
z^2$ on the unit circle $S^1$. As a lamination the dyadic solenoid
is one-dimensional; moreover,  it is naturally endowed with  the
structure of a compact topological group.

The map $f_0$ admits a natural extension $\h{f}_0$ acting on the
solenoid $\sol$. It turns out that $\h{f}_0$-periodic leaves in
$\sol$ are in one-to-one correspondence with periodic points of
$f_0$ in $S^1$. In turn, the periodic points of $f_0$ are in
one-to-one correspondence with rational numbers with odd
denominators modulo $1$. These observations are the first step
towards the proof of the Main Theorem.

To be more precise, let us introduce some objects that will play a
key role: We call a \textit{solenoidal cone} any space
homeomorphic to the inverse limit
$\cone=\inver(f_0,\h{\C}\setminus {\mathbb{D}})$. Due to
B\"ottcher's coordinate at infinity, for every quadratic
polynomial $f_c$ with connected Julia set, $Con_c$ the lift of the
basin of infinity to $\inver(f_c,\h{\C})$ is homeomorphic to
$\cone \setminus \sol$. We call $Con_c$ the \textit{solenoidal
cone at infinity} of $f_c$.

Now, note that $\cone$ is foliated by 1-dimensional solenoids
coming from the lifts of the equipotentials in the dynamical
plane. By means of B\"ottcher's coordinate, we can transfer this
foliation to $Con_c$. In $Con_c$ we call any of the corresponding
solenoids a \textit{solenoidal equipotential}. Each solenoidal
equipotential is canonically identified with $\sol$. Such
identification, allow us to keep track of periodic leaves in
$\mc{R}_c$ in a neighborhood of infinity.

A superattracting parameter $c$ admits a combinatorial model given
in terms of external rays, we will see that this combinatorial
model persists in the realm of the inverse limit of $f_c$. As a
matter of fact, this combinatorial model is the same for all
parameters inside a given hyperbolic component in the Mandelbrot
set. At the end of Section 3, we prove that the regular leaf space
associated to any hyperbolic quadratic polynomial $f_c$ is
homeomorphic to the regular leaf space of the center of the
hyperbolic component containing $c$. Hence, it is enough to
describe the topology of the centers of  hyperbolic components in
the Mandelbrot set.

\medskip

In Section 4, we describe the topology of the \textit{laminated
Julia set}, that is, the lift of the Julia set to the regular leaf
space. When $f_c$ is convex cocompact, the laminated Julia set is
compact, see \cite{LM}.  We show that if the postcritical set is
not the whole Julia set and the Julia set is locally connected,
then $f_c$ is convex cocompact if and only if the laminated Julia
set is leafwise connected.

By another result in \cite{LM}, if $f_c$ is convex cocompact, then
all the leaves of the regular leaf space $\mc{R}_c$ are
conformally isomorphic to the complex plane. Given a leaf
$L\subset{R}_c$, let us consider the set of \textit{unbounded
Fatou components} in $L$. We show that non-periodic leaves of
$\mc{R}_c$ have at most $2$ unbounded Fatou components. On the
other hand, the number of unbounded Fatou components on periodic
leaves depends on the valence of the corresponding repelling
periodic point.

Once we have described the basic topology of the laminated Julia
set, we are ready to give some restrictions on homeomorphisms
between regular leaf spaces of superattracting parameters. Because
the regular leaf space is locally compact and the laminated Julia
set is compact, we can compute the number of unbounded Fatou
components in some leaf $L$ in terms of the topology at infinity
of $L$.

Assume $f_c$ is superattracting, and let $U$ be the Fatou
component containing the critical point. The first return map of
$U$ has a unique fixed point, say $r_c$, on the boundary of $U$.
Following Milnor, $r_c$ is called the \textit{dynamic root} of
$f_c$. Accordingly, we call the cycle of $f_c$ containing the
dynamic root, the \textit{dynamic root cycle} of $f_c$. Since
$r_c\in J(f_c)$ and $f_c$ is superattracting, the point $r_c$ must
be a repelling periodic point of $f_c$. We prove that the topology
at infinity of the leaves containing the lift of the dynamic root
cycle is unique among all other leaves. Then, any homeomorphism as
in the Main Theorem must send leaves containing the lift of the
dynamic root cycle of $f_c$ into the corresponding leaves of
$\mc{R}_{c'}$.

\medskip

In Section 5, we prove the Main Theorem. The strategy of the proof
is to replace the homeomorphism $h:\mc{R}_c\ra \mc{R}_{c'}$ by
another homeomorphism $\tilde{h}$ with special characteristics.

First we show that $h$ is isotopic to a homeomorphism $\tilde{h}$
that sends a solenoidal equipotential $\mc{S}_{R}$ of $\mc{R}_c$
homeomorphically into a solenoidal equipotential $\mc{S}_{R'}$ of
$\mc{R}_{c'}$. Now, using the canonical identifications with
$\sol$ on $\mc{S}_R$ and $\mc{S}_{R'}$, the map $\tilde{h}$
induces a self-homeomorphism of the dyadic solenoid, which
according to Kwapisz \cite{Kwap} is isotopic to an \textit{affine
transformation} of the dyadic solenoid (that is, to the
composition of an automorphism and a translation). In order to get
$\tilde{h}$, we discuss some isotopic properties of
self-homeomorphisms of the dyadic solenoid $\sol$ and its cone
$\cone$.

The last part of the argument is to show that the translation part
of the affine transformation of the previous paragraph is isotopic
to zero. Thus after a series of improvements of the original
homeomorphism, we get a homeomorphism $\tilde{h}$ sending
$\mc{S}_R$ in $\mc{R}_c$ onto $\mc{S}_{R'}$, and such that the
restriction of $\tilde{h}$ to the solenoidal equipotentials is the
identity. It follows that the orbit portraits of the dynamic root
cycles of $c$ and $c'$ are the same, which completes the proof of
the Main Theorem.

There is an algebraic approach to laminations given by the theory
of Iterated Monodromy Groups developed by Volodymyr Nekrashevych
(see \cite{Nek} and \cite{BGN}). These Iterated Monodromy Groups
are precisely the monodromy groups of the regular spaces of
postcritically finite rational maps. For special cases, we were
able to prove some of the results in this paper using just
properties of the corresponding Iterated Monodromy Group (see
\cite{cab}). However, in general the combinatorics of these groups
is more subtle to describe. A good exposition of combinatorial
methods for rational maps can be found in Kevin Pilgrim's survey
\cite{Pil}.\\

\textbf{Acknowledgements.} I would like to thank Tomoki Kawahira
and Travis Waddington for helpful comments and discussions. I also
would like to thank the Fields Institute for their support and
hospitality while this paper was finished. This paper is based on
part of my thesis \cite{cab}. Finally, I would like to thank my
advisor Mikhail Lyubich.

\medskip

\section{Basic theory}

We assume the reader is familiar with the basic theory of dynamics
of rational maps. Definitions and properties of the objects
involved are widely spread in the literature; a very nice
exposition can be found in Milnor's book \cite{Mdyn}, and other
sources are \cite{CG} and \cite{L}.

We will review three combinatorial models of postcritically finite
quadratic polynomials, each one of which will be useful in this
work when we want to prove properties of regular spaces associated
to postcritically finite polynomials. For the proof of the Main
theorem we are going to use the model of ray portraits as
presented by Milnor in \cite{Mext}; for the interested reader
there is a complete discussion of combinatorics of postcritically
finite quadratic polynomials in \cite{BS}. At the end of the
section, we will discuss some properties for critically
non-recurrent quadratic polynomials.

\subsection{Combinatorics of postcritically finite polynomials}

A natural way to classify postcritically finite polynomials is by
describing the different arrangements that the critical orbit may
have as a subset in the plane. The combinatorics of this
arrangement are reflected in both dynamical and topological
properties of the Julia set. In fact, we could say that similar
properties hold for a slightly larger set of parameters; namely,
postcritically non-recurrent parameters. Combinatorics also make
explicit the relationship between the parameter and the dynamical
plane. So certain combinatorial behavior of a given parameter $c$
determines a region in the parameter plane where $c$ must lie.

\subsection{Hubbard trees}

To begin with, let us present a combinatorial model given by
Douady and Hubbard in \cite{DH}; for any postcritically finite
parameter $c$ they constructed an abstract graph, called the
\textit{Hubbard tree}, describing the dynamical arrangement of the
postcritical set. In \cite{DH}, Douady and Hubbard also proved
that different combinatorics induce different Hubbard trees.

The Hubbard tree is properly embedded in the dynamical plane as a
subset of the filled Julia set and can be defined by a finite set
of vertices in the Julia set. To see this, we note that between
any two points $z$ and $\zeta$ in the Julia set $J(f_c)$ there is
an unique arc $\gamma$, embedded in the Julia set, connecting $z$
with $\zeta$. The uniqueness of $\gamma$ is subject to the
condition that, if the trajectory of $\gamma$ intersects a Fatou
component, then it goes along internal rays. In this way, the
\textit{Hubbard tree} of a postcritically finite quadratic
polynomial $f_c$ is defined as the smallest collection of arcs,
embedded in the Julia set $J(f_c)$ and connecting the entire
critical orbit. This tree is finite and forward invariant under
the action of $f_c$.

\subsection{Orbit portraits}

The description by orbit portraits of postcritically finite maps
is based upon the following proposition:

\begin{proposition}[Douady and Hubbard]\label{dyn.land.ray}
Every repelling and parabolic periodic point of a quadratic
polynomial $f_c$ is the landing point of an external ray with
rational angle. Conversely, every external ray with rational angle
lands either at a periodic or preperiodic point in $J(f_c)$.
\end{proposition}

From now on, we will follow the exposition of John Milnor in
\cite{Mext}, see also the related work of Alfredo Poirier
\cite{Poi} and Dierk Schleicher \cite{Schext}. By Proposition
\ref{dyn.land.ray}, every parabolic or repelling cycle
$P=\{p_1,...,p_n\}$ is associated to a family of finite sets
$\mc{O}_P=\{A_1,...,A_n\}$, where $A_i=\{\theta\in
\mathbb{Q}/\mathbb{Z}| R_\theta$ lands in $p_i\}$. The family
$\mc{O}_P$ is called the \textit{orbit portrait} of the cycle $P$.
The doubling map $\mc{D}$ on $\mathbb{T}$ permutes the sets $A_i$
and acts on the angles in $A_i$ in an order preserving way.
Moreover, every angle in $\cup A_i$ is periodic under $\mc{D}$
and the period of such angle only depends on the cycle. As a
consequence, all $A_i$ have the same cardinality. Thus, for $p\in
P$, the \textit{valence} $v_p$ of $p$ is the cardinality of any of
the sets $A_i\in \mc{O}_P$.

Given two angles $\theta_1$ and $\theta_2$ in the unit circle
$\mathbb{T}$, by $\widehat{\theta_1\theta_2}$ we denote the
\textit{directed arc} in $\mathbb{T}$ from $\theta_1$ to
$\theta_2$. Now, if $A$ is a finite set in $\mathbb{T}$, a
\textit{complementary arc} of $A$ is the closure of any connected
component of $\mathbb{T}\setminus A$. If the valence is bigger
than one, each $A_i$ determines a finite collection of
complementary arcs in the circle. Among the union of all
complementary arcs of all $A_i$, see Lemmas 2.5 and 2.6 in
\cite{Mext}, there is a unique arc of shortest length called the
\textit{characteristic arc}. The preimage of the characteristic
arc under doubling is also the longest complementary arc which is
called the \textit{critical arc}; it is worth noting that the
critical arc can also be defined as the unique complementary arc
subtending an angle bigger than $1/2$. Since the doubling map
preserves the order in the complementary arcs, if
$\widehat{\theta_1 \theta_2}$ is the characteristic arc of some
orbit portrait $\mc{O}_P$, then by a straight forward calculation,
the critical arc is $\widehat{\eta_1\eta_2}$ where
$\eta_1=\theta_1/2+1/2$ and $\eta_2=\theta_2/2$.

\subsection{The dynamic root point}

Let $c$ be a superattracting parameter of period $n$, then the
first return map of the Fatou component containing $c$ has a
unique fixed point on the boundary. Following Milnor \cite{Mext}
and Schleicher \cite{Schext}, we call $r_c$ the \textit{dynamic
root} of the superattracting parameter $c$. The orbit portrait of
$r_c$ is called the \textit{critical portrait} of $c$. It is a
Theorem by Milnor \cite{Mext} that, the critical portrait
characterizes the parameter $c$. In other words, no two
superattracting parameters have the same critical portrait (see
also \cite{Poi}). Recall that rotations on the unit circle
$\mathbb{T}$ are given by maps $r_{\theta}:\mathbb{T}\rightarrow
\mathbb{T}$ of the form $r_{\theta}(\tau)=\tau +\theta$, mod
($\mathbb{Z}$). The following lemma will be important for the
closing argument in the Main Theorem:

\begin{lemma}\label{orbit} Let $\mc{O}_P$ and $\mc{O}_{P'}$ be
the orbit portraits of the periodic cycles $P$ and $P'$. If there
is a rotation of the circle $r_\theta$, such that
$r_\theta(\mc{O}_P)=\mc{O}_{P'}$, then $\mc{O}_{P}=\mc{O}_{P'}$.
In particular, if $\mc{O}_{c}$ and $\mc{O}_{c'}$ are the critical
portraits of two superattracting quadratic polynomials differing
by a rotation, then $c=c'$
\end{lemma}

\begin{proof} Since the characteristic arc is the minimal complementary
arc, the rotation $r_{\theta}$ must send the characteristic arc of
$P$ to the characteristic arc of $P'$. Analogously, $r_{\theta}$
must send the critical arc of $P$ to the critical arc of $P'$. If
$\widehat{\theta_1 \theta_2}$ and $\widehat{\theta'_1 \theta'_2}$
are the characteristic arcs of $P$ and $P'$ respectively, then the
hypothesis yields the following equations:
$$\theta_i+\theta=\theta'_i$$ for the characteristic arcs, and
$$\theta_i/2+1/2+\theta=\theta'_i/2+1/2$$ for the critical arcs. Thus,
$\theta=0$ and  $\theta_i=\theta'_i$ for $i=1,2$, then the
critical arc of $f_c$ is equal to the critical arc of $f_{c'}$.
Since the characteristic arcs generate the orbit portrait, $P=P'$.
The second part of the lemma follows from the fact that the
critical portrait determines the parameter $c$.
\end{proof}

\subsection{Combinatorics in the parameter plane}

Now, we briefly discuss how the location of a certain parameter
$c$ in the Mandelbrot set affects the combinatorics of periodic
orbits in the dynamical plane of $f_c$. Analogous to Proposition
\ref{dyn.land.ray}, we have:

\begin{proposition}[Douady and Hubbard]\label{Mandelbrot.rays}
In the parameter plane, every parabolic or Misiurewicz parameter
is the landing point of at least one external ray of rational
angle. Inversely, every external ray with rational angle $\theta$
lands at some point $c$ in the boundary of the Mandelbrot set.
Moreover, if $\theta$ has odd denominator, then $c$ is a parabolic
parameter and Misiurewicz otherwise.
\end{proposition}

In this way, every rational angle is associated to either a
parabolic or a Misiurewicz parameter. Besides $c=1/4$, which
corresponds to the cusp of the main cardioid, every parabolic
parameter $c$ is the landing point of exactly two external rays,
say $R_{\theta_1}$ and $R_{\theta_2}$. These external rays cut the
plane into two parts, one of them contains the Main Cardioid and
the other is called the \textit{wake} $W_c$ determined by $c$. The
\textit{limb at} $c$ is the intersection $W_c\cap M$.

Every parabolic parameter $c$ can be regarded as the root of some
hyperbolic component $H_c$; if $c\neq 1/4$, then $H_c$ is
contained in $W_c$. Thus, except for the hyperbolic component
inside the Main Cardioid, the root of every hyperbolic component
$H$ disconnects $H$ from the Main Cardioid. Even more, and here is
the one of the most beautiful features of the combinatorics of
quadratic polynomials, if $R_{\theta_1}$ and $R_{\theta_2}$ are
the rays determining the wakes of the parabolic parameter $c$,
then $\widehat{\theta_1\theta_2}$ corresponds to the
characteristic arc of the parabolic cycle of $c$. Moreover,
$\widehat{\theta_1\theta_2}$ is also the characteristic arc of the
dynamic root cycle of $H_c(0)$, the center of the hyperbolic
component $H_c$. As for the whole wake $W_c$, it can be described
as the set of parameters in $M$ for which there is a cycle $P$
with orbit portrait $\mc{O}_P$, and such that the characteristic
arc of $\mc{O}_P$ is $\widehat{\theta_1 \theta_2}$. To any given
hyperbolic component $H$, we associate to $H$ the smallest angle
of the external rays landing at the root of $H$. Another useful
fact is that the period of $\theta_1$ is the same as the period of
the parabolic cycle of $c$ and the critical cycle of $H_c(0)$.

When $c$ is Misiurewicz either $c$ is the landing point of several
external rays or exactly one external ray. On the other hand, the
set of Misiurewicz parameters where exactly one ray lands is
countable.

\subsection{Internal address}

We now follow the discussion of Eike Lau and Dierk Schleicher in
\cite{LS}  to introduce internal addresses. Although the
definition of internal addresses is for hyperbolic components, we
can extend this definition to centers of the corresponding
parameters.

Given two hyperbolic components or, for this discussion,
Misiurewicz points, $H_1$ and $H_2$, we say that $H_2$ is
\textit{visible} from $H_1$ if the rays landing at the root of
$H_1$ separate $H_2$ from the Main Cardioid. Visibility induces a
partial ordering on the hyperbolic components. More precisely, for
every pair of hyperbolic components, or Misiurewicz parameters,
$H_1$ and $H_2$ either one is visible from the other, or there is
a hyperbolic component, or Misiurewicz parameter, $H_3$ from which
both $H_1$ and $H_2$ are visible.

Let $H$ be a hyperbolic component, then using visibility between
hyperbolic components and Misiurewicz parameters, it is possible
to define a path, called \textit{the combinatorial arc} of $H$,
along the Mandelbrot set connecting the Main Cardioid with $H$.
This imposes a tree structure onto the arrangement of hyperbolic
components in the Mandelbrot set. The root of this tree is the
Main Cardioid, and the set of vertices is the set of hyperbolic
components and Misiurewicz parameters.

In general, the combinatorial arc of $H_c$ crosses infinitely many
hyperbolic components. Nevertheless, one can write down, in an
increasing sequence of numbers $\{a_n\}$ the periods of the
hyperbolic components that the combinatorial path of $H_c$
crosses. It turns out, that this sequence $\{a_n\}$ is always
finite and is called the \textit{internal address} of the
component $H_c$, or the associated parameter $c$. Following Lau
and Schleicher's notation, we consider the component inside the
Main Cardioid as part of the combinatorial arc, so the internal
address always starts with $1$, and when we write down the
sequence of numbers $a_n$, we connect consecutive numbers by
arrows.

For example, the internal address of the ``\textit{basillica}" map
$f_{-1}=z^2-1$ is $1\rightarrow 2$, while the ``\textit{airplane}"
map $f_c(z)=z^2 +c$ with parameter $c=-1.7548776662466927601$ has
internal address given by $1\rightarrow 2 \rightarrow 3$. It is a
Theorem by Schleicher that if two hyperbolic parameters $c$ and
$c'$ have the same internal address then the maps restricted to
the filled Julia sets are conjugate. So, to make internal address
a efficient combinatorial model it is necessary to include
additional information. Namely, the rotation number around
periodic cycles:

Let $c$ be a parabolic parameter different from $1/4$, then the
multiplier of the corresponding parabolic cycle is a root of unity
of the form $e^{\frac{p}{q} \pi i}$, the number $p/q$ is called
the \textit{combinatorial rotation number} of $c$.

If, in the internal address, we label each arrow with the
combinatorial rotation number of the hyperbolic components on the
sequence, then we obtain the \textit{labelled internal address},
is a Theorem of Lau and Schleicher \cite{LS}, that the labelled
internal address of a superattracting parameter $c$ characterizes
the parameter.

\subsection{Postcritically non-recurrent quadratic polynomials}

Besides superattracting, hyperbolic and Misiurewicz parameters,
there is also a class of parameters such that the orbit of the
critical point does not contain the critical point in its
accumulation set. More precisely, let $f:X\rightarrow X$ be a
dynamical system defined in a metric space $X$. Let us define the
omega limit $\omega(x)$\textit{-limit set} of a point $x$ as the
set of accumulation points of $\{f^n(x)\}_{n\in \mathbb{N}}$ in
$X$. A non-periodic point $x_0 \in X$ is said to be a
\textit{recurrent} point of $f$ if $x_0 \in \omega(x_0)$. The
action of $f$ on a set $A$ is said to be \textit{non-recurrent} if
no point $a$ in $A$ is a recurrent point of $f$.

A stronger notion of recurrence in a set $A$ is when every point
of $A$ is an accumulation point of its orbit: Let $f:X\rightarrow
X$ be a dynamical system defined in a metric space $X$. A set $A$
is called \textit{minimal} if it is closed invariant under $f$ and
no proper subset of $A$ has this property.

Note that every point in a minimal set most be recurrent. A
Theorem of Birkhoff shows that every dynamical system contains a
minimal set.

A \textit{postcritically non-recurrent} rational map $f$ is a map
whose action on the postcritical set is non-recurrent. In the case
of quadratic polynomials, it has been proved by Carleson, Jones
and Yoccoz \cite{CJY}, that any postcritically non-recurrent
quadratic polynomial has locally connected Julia set. Another
non-trivial Theorem by Yoccoz states that critically non-recurrent
parameters are locally connected in the Mandelbrot set. In
\cite{CJY} these parameters are called \textit{subhyperbolic};
however we will introduce later another term related to the
associated lamination.

A Theorem by Fatou states that if all the critical points of a
rational map $f$ belong to $F(f)$ then, for every $x$ in $J(f)$
there is a number $C>0$ and  $\sigma>1$ such that $|(f^n)'(x)|>C
\sigma^n$. In other words, the map $f$ is expanding on the Julia
set. This is the case for all quadratic polynomials with
hyperbolic parameter. Now, the following Theorem by Ricardo
Ma\~{n}e \cite{Man} describes, in the general case, those points
in the Julia set which are expanding.

\begin{theorem}[Ma\~{n}e's Theorem] Let $f:\bar{\C} \rightarrow
\bar{\C}$ be a rational map. A point $z\in J(f)$ is either a
parabolic periodic point, or belongs to the $\omega$-limit set of
a recurrent critical point, or for every $\epsilon >0$, there
exists a neighborhood $U$ of $x$, such that, $\forall n \geq 0$
every connected component of $f^{-n}(U)$ has diameter $\leq
\epsilon$.
\end{theorem}

A related, and useful, result is the following Lemma:

\begin{lemma}[Shrinking Lemma] \label{shrinking} Let $f: \bar{\C} \rightarrow \bar{\C}$ be a
rational map. If $K\subset J(f)$ is a compact subset disjoint from
parabolic periodic points and $\omega$-limit sets of recurrent
critical points, then for every $\epsilon>0$ there exist
$\delta>0$ such that for every $k\in K$ and every $n\geq 0$, all
connected components of $f^{-n}(B(k,\delta))$ have diameter less
than $\epsilon$.
\end{lemma}

The wilder the postcritical set of a function $f$ is, the more
intricate becomes the description of the dynamics of $f$. In this
work, when speaking of critically recurrent quadratic polynomials
$f_c$, we will consider only parameters with locally connected
Julia set, and such that the postcritical set $P_c$ is a Cantor
set and the action of the map $f_c$ restricted to $P_c$ is
minimal.



\section{Inverse limits and laminations}

In this work we follow the notation and definitions as found in
the papers \cite{LM}, \cite{KL} and \cite{L}. A topological space
$\mc{B}$ is said to have a \textit{product structure} if it is
provided with a homeomorphism $\phi:\mc{B}\rightarrow \mathbb{D}^n
\times T$ where $\mathbb{D}^n$ denotes the open unit disk in
$\mathbb{R}^n$, and $T$ is some topological space. Sets of the
form $B_x=\phi^{-1}(\mathbb{D}^n\times \{x\})$ are called
\textit{local leaves} or simply \textit{plaques}, while sets of
the form $T_z=\phi^{-1}(\{z\}\times T)$ are called \textit{local
transversals}.

A \textit{lamination} is a Haussdorf topological space $\mc{X}$
which is endowed with an atlas of open charts $(\phi, U)$ where
$U$ has a product structure and $\phi$ is a homeomorphism as
above.  We also require that the change of coordinates are laminar
maps, thus change of coordinates are maps of the form $\gamma
_{\alpha \beta }:D\times T\rightarrow D^{\prime }\times T^{\prime
}$ given by $\gamma_{\alpha \beta} (z,t)=(\sigma (z,t),\psi (t)),$
where $\sigma $ and $\psi$ are continuous functions on $t$. The
sets $U$ will be called \textit{flow boxes}. A \textit{laminar
map} between laminations is a continuous map such that when
restricted to a flow box it sends plaques into plaques.

Different regularity conditions can be imposed on laminations.
This is done by requiring the corresponding regularity of the map
$\sigma$ along $z$. For instance, smooth laminations are
laminations whose transition maps $\gamma_{\alpha,\beta}$ are
smooth in the $z$ variable. Similarly, real and complex analytic
laminations can be defined. Laminations are generalizations of the
concept of foliations; a foliation is a lamination where $\mc{X}$
is a manifold itself.

A lamination $\mc{X}$ is decomposed into a disjoint union of
connected $n$-manifolds $\sqcup L_{\alpha }$, where the sets
$L_{\alpha}$ are called \textit{global leaves}, or just
\textit{leaves}. A global leaf can be characterized as the
smallest set $L$ with the property that if it intersects a plaque
$B_x$ then $B_x\subset L$. The maps $\phi$ restricted to plaques
are, in fact, charts for the leaves. Given a point $z$ in
$\mc{X}$, we will denote by $L(z)$ the leaf in $\mc{X}$ which
contains $z$.

We define the \textit{dimension of a lamination} $\mc{X}$ as the
dimension of any plaque in $\mc{X}$. In dimension two, the concept
of a conformal lamination is equivalent to the one of a complex
analytic lamination, these type of laminations are also called
\textit{Riemann surfaces laminations}.

Other categories of laminations, which play an important role, are
affine and hyperbolic laminations. For these, we require the
change of coordinates on leaves to be affine and hyperbolic
isometries, respectively. In this work, we will be interested in
the topology of affine laminations of dimension two arising from
dynamics of quadratic polynomials.

\subsection{Inverse limits}

Consider $\{f_k:X_k \rightarrow X_{k-1}\}$, a sequence of
$m$-to-$1$ branched covering maps between $n$-manifolds $X_k$.
Then, define the \textit{inverse limit} $\inver(f_n,X_n)$ as
$$\inver (f_n,X_n)=\{\h{x}=(x_1,x_2,...)\in \prod X_n |f_{n+1}(x_{n+1})=x_n\}.$$ The space $\inver(f_n,X_n)$ has a \textit{natural topology} which
is induced from the product topology in $\prod X_n$.

We are interested in inverse limits arising from dynamics; these
are particular cases where all coverings $f_n \equiv f$ and the
manifolds $X_n$ are equivalent to a single phase space $X$. Such
inverse limits are called \textit{solenoids}, or \textit{natural
extensions} in ergodic theory, we will denote them by
$\inver(f,X)$ to make emphasis in $X$. When $f$ is a rational
function and $X=\bar{\C}$ then, following Lyubich and Minsky
\cite{LM}, we will denote $\inver(f,\bar{\C})$ by $\mc{N}_f$.

The map $f:X\rightarrow X$ has a \textit{natural extension}
$\h{f}:\inver(f,X) \rightarrow \inver(f,X)$ defined as
$$\h{f}(x_0,x_{-1},...)=(f(x_0),x_0,x_{-1},...).$$ Also, there is a
family of natural projections $\pi_{-n}:\inver(f,X) \rightarrow
X$, given by $\pi_{-n}(\h{x})=x_{-n}.$ Each of these maps
semiconjugates $\h{f}$ to $f$, so $\pi_{-n}(\h{f}(\h{z}))=
f(\pi_{-n}(\h{z}))$. For simplicity, the subindex of the
projection over the first coordinate will be omitted, thus
$\pi\equiv \pi_0$. We are interested in studying properties of
dynamics of $\h{f}$ and how they are related to the dynamics of
the original map $f$.

Let $A$ be an invariant set of $f$, so we have $f(A)\subset A$.
The \textit{invariant lift} of $A$ in $\inver(f,\bar{\C})$ is the
set $\h{A}$ of all backward orbits $\h{z}$ such that
$\pi_{-n}(\h{z})\in A$ for every $n$. Clearly invariant lifts of
periodic cycles of $f$ are periodic cycles of $\h{f}$.

Note that there is a natural one-to-one identification of the
periodic points of $f$ with the periodic points of $\h{f}$ in
$\mc{N}_f$. Namely, to every periodic point $p$ of $f$, the
corresponding point $\h{p}$ is the point in the invariant lift of
the cycle of $p$ such that $\pi(\h{p})=p$. Vice versa, given a
periodic point $\h{p}$ of $\h{f}$, the point $\pi(\h{p})$ is a
periodic point of $f$. When $f$ is a rational function, we
classify invariant lifts of periodic cycles by borrowing the
corresponding classification in the dynamical plane. For instance,
a periodic point $\h{p}$ in $\inver(f,\bar{\C})$ is called
parabolic if $\pi(\h{p})$ is parabolic in $\C$.

\subsection{Lamination structure of solenoids}

Assume $f$ does not have critical points, so it is an $m$-to-$1$
covering map, and let $\h{x}=(x_0,x_{-1},x_{-2},...)$ be a point
in $\inver(f,X)$;  each coordinate $x_n$ belongs to the set of $m$
preimages of $x_{n-1}$. Hence, after a suitable labelling of the
branches of $f^{-1}$, the fiber $\fiber{(x_0)}$ can be identified
with $\{0,...,m-1\}^ {\mathbb{N}}$; moreover, by taking the
product discrete topology on $\{0,...,m-1\}^ {\mathbb{N}}$, this
identification is a homeomorphism. Let $U$ be an open neighborhood
of $x_0$. The set $\fiber{(U)}$ is an open neighborhood around
$\tilde{x}$ and is homeomorphic to $U\times \{0,...,m-1\}^
{\mathbb{N}}$. This endows $\inver(f,X)$ with a lamination
structure. For $f$ with critical points, whenever an open set
$U\subset\C$ does not intersect the postcritical set of $f$, the
set $\pi^{-1}(U)$ has a product structure.

Let $(U_{-n})$ be the \textit{pull-back of} $U=U_0$ \textit{along}
$\h{x}$, where $U_{-n}$ is the connected component of $f^{-n}(U)$
containing $x_{-n}$. Given a number $N$, the sets
$$B(U,\h{z},N)=\pi^{-1}_{-N}(U_{-N})$$ form a local basis of open sets
for $\inver(f,X)$. A plaque, then, can be regarded either as a
connected component of a flow box, or as the complete pull back of
some open disk $U_0$ along $\h{x}$; that is, a sequence of the
form $(U_0,U_{-1},U_{-2},...)$.

 \subsection{Dyadic solenoid}

\begin{figure}[htbp]
  \begin{center}
  \includegraphics[width=2 in]{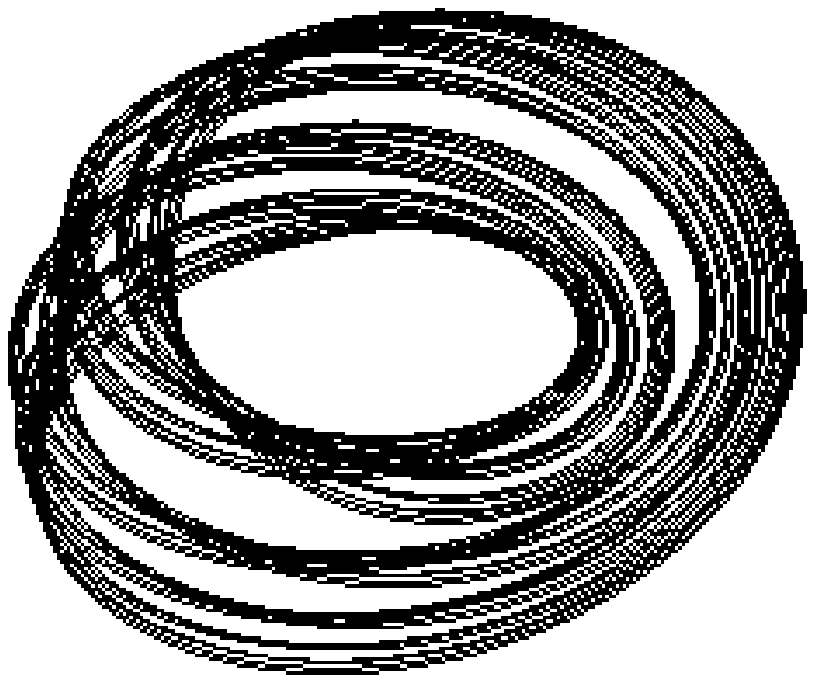}
  \caption{The dyadic solenoid $\sol$.}
  \label{fig.dyad.solenoid}
  \end{center}
\end{figure}

Consider the polynomial $f_0(z)=z^2$ defined on $X=\mathbb{S}^1$.
The solenoid $\mc{S}^1=\inver(f_0,\mathbb{S}^1)$ is called the
\textit{dyadic solenoid}, see Figure~\ref{fig.dyad.solenoid}. As
$f_0$ is a covering map of degree two, for every point $z\in
\mathcal{S}^1$, the fiber $\pi^{-1}(\pi(z))$ is homeomorphic to
$\cantor$. Since $\mathbb{S}^1$ is a compact topological group,
the solenoid $\sol$ is a compact topological group, in which the
group multiplication is given by the multiplication of
$\mathbb{S}^1$ applied componentwise, so the unit $\h{u}$ is the
point $(1,1,1,...)$. The translations of the solenoid, given by
left multiplication of elements in $\sol$, will be denoted by
$\tau_{\h{\zeta}}(\h{z})=\h{\zeta}\cdot \h{z}$.

The leaf containing the unit $\h{u}$ in $\sol$ is a one-parameter
subgroup of the dyadic solenoid, parameterized by the map
$\rho:\mathbb{R}\rightarrow \sol$, with $\rho(t)=(e^{2\pi i
t},e^{\pi i t},e^{\pi i t/2},...)$. The image of $\rho$ is dense
in the solenoid, that is $\sol=\overline{\rho(\mathbb{R})}$, and
since $\sol$ is a topological group every leaf in the dyadic
solenoid is dense. Since in fact, there are uncountable many
leaves in the solenoid, let us remark that the dyadic solenoid is
connected but not path connected. By transferring the natural
order in $\mathbb{R}$ to the leaves in $\sol$, the map $\rho$ also
introduces a leafwise order in $\sol$, namely, if $\h{z}$ and
$\h{\zeta}$ are two points in the same leaf, then we say that
$z>\zeta$ whenever $\rho^{-1}(z\cdot\zeta^{-1})>0$.

\subsubsection{The doubling map}

Let us consider the polynomial $f_0=z^2$ which is associated to
the center of the component inside the Main Cardioid. The orbit of
any point inside the unit circle converges to $0$ under iteration
of $f_0$; while the orbit of any point outside the unit circle
tends to $\infty$. Hence, the Fatou set consists of only two
domains and, the Julia set is just the unit circle $\mathbb{S}^1$
in the complex plane $\C$.

The action of $f_0$ in the unit circle is $f_0(e^{2\pi
i\theta})=e^{2(2\pi i \theta)}$. That is, it doubles the
corresponding angle $\theta$. So, the map
$f_0:\mathbb{S}^1\rightarrow \mathbb{S}^1$ is conjugate to the
doubling map $\mc{D}:\mathbb{R}/\mathbb{Z}\rightarrow
\mathbb{R}/\mathbb{Z}$, which sends the angle $\theta$ to
$2\theta$ (mod $1$). For convenience, whenever we refer to a point
in the unit circle, we denote it by its angle in
$\mathbb{T}=\mathbb{R}/\mathbb{Z}$. We also assume the standard
orientation on the unit circle. The orbit under doubling of every
rational angle is finite, and the set of periodic points of
$\mc{D}$ is exactly the set of rational angles $p/q$ where $q$ is
an odd number.

Using the argument above, it is easy to check that the dyadic
solenoid $\sol$ is isomorphic to the solenoid
$\inver(\mc{D},\mathbb{T})$. Thus, the periodic points of
$\h{f_0}$ are in one-to-one correspondence with the reduced
rational angles $p/q$ in $\mathbb{Q}/\mathbb{Z}$ where $q$ is an
odd number.

By parameterizing external rays with their corresponding angle in
$\mathbb{T}$, and as a consequence of B\"{o}ttcher conjugacy of
$f_c$ to $z^2$ on the basin of infinity, we have
$f_c(R_\theta)=R_{2\theta}$. So, the action of $f_c$ on the set of
external rays is conjugate to the action of the doubling map in
$\mathbb{T}$. This feature will play a very important role in this
work.

\subsubsection{The adding machine}

The projection $\pi:\inver(f_0,\mathbb{S}^1)\rightarrow
\mathbb{S}^1$ is a fibration map, the monodromy group of this
fibration is isomorphic to $\mathbb{Z}$. As we mentioned above,
the fiber of every point in $\mathbb{S}^1$ is homeomorphic to the
set $\cantor$ which, in turn, can be identified with the set of
formal series of the form $\sum^{\infty}_{i=0} a_i 2^i$ with
$a_i\in \{0,1\}$.  Given such identification, the action of the
Monodromy group is generated by the \textit{adding machine} on
$\cantor$, this is the map given by adding one to an element in
$\cantor$. Let us take $F$ to be the fiber of $\pi$ over $1$. The
adding machine action on $F$ extends continuously to the whole
solenoid. In fact, its generator is given explicitly by the map
$\sigma:\sol \rightarrow \sol$, defined by
$\sigma(\h{z})=\rho(1)\cdot \h{z}$. We refer to the action of the
group $<\sigma>$ as the \textit{adding machine action}. Finally,
let us note that the solenoid $\sol$ is the quotient of the
product space $S=[0,1]\times F$ by the relation
$(1,f)\sim(0,\sigma(f))$.

Two points $\h{z}$  and $\h{\zeta}$ in $F$ belong to the same leaf
if and only if $\h{z}$ and $\h{\zeta}$ belong to the same orbit
under the adding machine action. As $F$ is homeomorphic to
$\cantor$, $F$ is uncountable, but each orbit under the adding
machine action has a countable number of points in $F$, so the
dyadic solenoid has uncountable many leaves.

\subsection{Solenoidal cones}

\begin{figure}[tbph]
  \begin{center}
  \includegraphics[width=2 in]{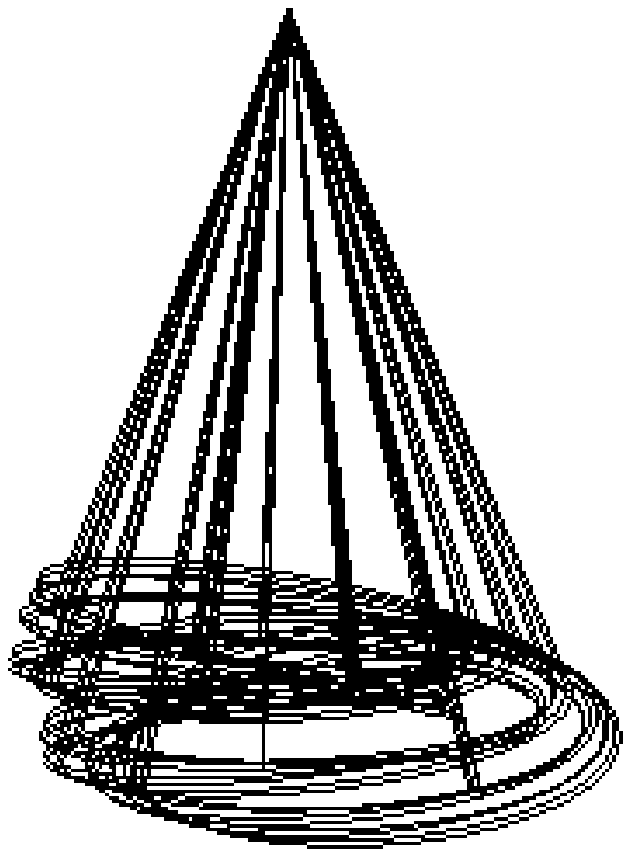}
  \caption{The solenoidal cone $\cone$.}
\label{solenoidalcone}
  \end{center}
\end{figure}

The solenoids $\inver(f_0,\mathbb{D}^*)$ and
$\inver(f_0,\mathbb{C}\setminus \bar{\mathbb{D}})$, are both
homeomorphic to $\mc{S}^1\times (0,1)$. Because $\infty$ is a
superattracting fixed point of $f_0$, the inverse limit
$\inver(f_0, \bar{\C} \setminus \mathbb{D})$ is homeomorphic to
the cone over the dyadic solenoid $\mc{S}^1$ defined as
$\sol\times[0,1]/\{(s,1)\sim(s',1)\}$ for all $s,s'\in \sol$. The
vertex of $\inver(f_0, \bar{\C} \setminus \mathbb{D})$ corresponds
to the point $\h{\infty}=(\infty,\infty,\infty,....)$, we call any
space homeomorphic to $\inver(f_0, \bar{\C} \setminus \mathbb{D})$
a \textit{closed solenoidal cone}. In particular, the solenoidal
cone $\inver(f_0,\bar{\C}\setminus \mathbb{D})$  will be denoted
by $\cone$. The dyadic solenoid $\sol$ is contained in $\cone$,
and we regard $\sol$ as the boundary of $\cone$ (see
figure~\ref{solenoidalcone}). Since there is no local product
structure on $\h{\infty}$, the solenoid $\cone$ is not a
lamination. In general, for dynamical systems with critical
points, the situation where critical points occur infinitely many
times in the coordinates of a given point is one of the possible
obstructions for the inverse limit to be a lamination.

It is important to note that $\inver(f_0,\C \setminus
\mathbb{D})\simeq \sol \times (0,1)$ is not path connected,
whereas $\cone$ is path connected.

As a consequence of B\"ottcher's Theorem, we can associate to
every quadratic polynomial $f_c=z^2+c$ a solenoidal cone. When $c$
belongs to the Mandelbrot set, the B\"{o}ttcher's coordinate
$\phi_c:(A(\infty),\infty) \rightarrow (\bar{\C}\setminus
\bar{\mathbb{D}},\infty)$ conjugates $f_c$ to $f_0$. The map
$\phi_c$ naturally lifts to a homeomorphism
$\h{\phi}_c:\inver(f_c,A(\infty))\rightarrow \cone \setminus \sol$
given by
$\h{\phi}_c(z_0,z_{-1},...)=(\phi_c(z_0),\phi_c(z_{-1}),...)$
which conjugates $\h{f}_c$ to $\h{f}_0$. We will call the solenoid
$\inver(f_c,A(\infty))\subset \mc{N}_{f_c}$ the \textit{solenoidal
cone at infinity} of $f_c$. When $c$ does not belong to the
Mandelbrot set, we can associate a solenoidal cone at infinity to
$f_c$, however this cone does not project to the whole basin of
infinity but to a neighborhood of $\infty$. We will construct such
a cone in next subsection.

The inverse limit $\inver(f_0,\bar{\mathbb{D}})$ is also a
solenoidal cone, both $\cone$ and $\inver(f_0,\bar{\mathbb{D}})$
share the common boundary $\sol$. Thus, the natural extension
$\mc{N}_{f_0}$ can be decomposed into $\cone$ and
$\inver(f_0,\bar{\mathbb{D}})$ by cutting along $\sol$. This shows
that $\mc{N}_{f_0}$ is homeomorphic to the double cone over the
solenoid $\mc{S}^1$.

Note that the one-dimensional solenoid $\sol$ is connected and
every leaf in it is dense. Thus, the solenoids $\sol\times (0,1)$
and $\cone$ are connected, and every leaf in each of them is
dense.

\subsection{Subsolenoidal cones}

Let $r>1$, and $\mathbb{D}_r=\{z\in \C| |z|<r\}$, then there is a
\textit{canonical homeomorphism} between $\inver(f_0,\C\setminus
\mathbb{D}_r)$ and $\sol\times [0,1)$ given by $$(z_0,
z_{-1},...)\mapsto (\frac{z_0}{|z_0|},\frac{z_{-1}}{|z_{-1}|},...
)\times (1-\frac{r}{|z_0|}),$$ which extends to a homeomorphism
between $\inver(f_0,\bar{\C}\setminus \mathbb{D}_r)$ and $\cone$,
fixing $\h{\infty}$. This implies that, for $f_c$ and the space
$A_r=\{z\in A_c(\infty)| |\phi_c(z)|\geq r\}$ outside the
equipotential $E_r$ in $A_c(\infty)$,  the solenoid
$\inver(f_c,A_r)$ is also homeomorphic to $\cone$.

As we mention above, when $c$ does not belong to the Mandelbrot
set, we can not associate a solenoidal cone to the lift of the
basin of infinity, however for $r$ big enough, we can take
$\inver(f_c,A_r)$ as the solenoidal cone at infinity for $f_c$.

The fiber over $\pi^{-1}(E_r)$ of any equipotential $E_r$ is
homeomorphic to the dyadic solenoid $\mc{S}^1$. We call such a
space, the \textit{solenoidal equipotential} over $E_r$, and
denote it by $\mc{S}_r$. Solenoidal equipotentials form a
foliation of the solenoidal cone at infinity in $\mc{N}_{f_c}$.
Similarly, external rays induce a foliation in
$\pi^{-1}(A(\infty))$ of \textit{solenoidal external rays}, or
just \textit{solenoidal rays} for short. The fiber over each
external ray is a ``Cantor set" of external rays.

In particular, in the fiber of a periodic external ray under
$f_c$, exactly one of the components is periodic under $\h{f}_c$.
Moreover, every periodic solenoidal ray must project to a periodic
external ray in the dynamical plane. Hence, periodic solenoidal
external rays are parameterized by reduced rational numbers
$\frac{p}{q}\in \mathbb{T}$ with $q$ odd. Thus, to every periodic
solenoidal ray we can associate a periodic point in $\sol$. It
turns out, that after identifying a solenoidal equipotential
$\mc{S}_r$ with $\sol$, each periodic solenoidal rays intersect
any solenoidal equipotential $\mc{S}_r$ at a periodic point in
$\sol$.

Moreover, every periodic solenoidal ray lands at a periodic point
$\h{p}$ of $\h{f}_c$ in $\mc{N}_c$. So if two external rays land
at $\pi(\h{p})$ then the corresponding periodic lifts of these
rays belong to the same leaf. In this way, external ray
identification corresponds to an identification of different
leaves in $\sol$.

%

\subsection{Lyubich-Minsky laminations}

Finally, we get to the object of study of this paper. Inverse
limits are also defined for branching mappings, however, as we saw
in $\cone$, these solenoids may not have a local lamination
structure. The case for rational maps on the sphere was addressed
by Lyubich and Minsky in \cite{LM}. In this setting, there are two
obstructions for the natural extension to be a lamination. The
first, as we have noted before, is structural; there are some
points which fail to have local product structure. The second
comes from the fact that the geometric structure on the leaves may
not vary continuously; on the transverse direction. Thus, it is
necessary to refine the topology in an appropriate way. The first
problem is easy to carry over, simply by removing all points that
don't have local product charts. The second is of a more delicate
nature; in fact, for many rational functions the new topology is
hard to describe by intrinsic properties of the natural extension.
The way out found by Lyubich and Minsky was to embed the
lamination into a universal space where a suitable topology comes
naturally.

First let us introduce the regular space. To keep things simple,
we return to the assumption that $f_c$ is a quadratic polynomial
defined on the Riemann sphere.

\begin{definition} A point $\h{z} \in \inver(\bar{\C},f_c)$ is called
\textit{regular} if there exists a neighborhood $U_0$ of
$\pi(\h{z})$ such that the pull back of $U_0$ along $\h{z}$ is
eventually univalent. The set $\mc{R}_{c}\subset
\inver(\bar{\C},f_c)$ of regular points is called the
\textit{regular part} of $f_c$.
\end{definition}

The regular part is a disjoint union of Riemann surfaces. To check
this, note that, in terms of plaques
$\h{U}=(U_0,U_{-1},U_{-2},...)$, a point $\h{z}$ is regular if
there is a number $N\geq 0$ such that $U_{-n}$ does not contain
critical points for $n\geq N$. Then, naturally, for a conformal
chart for the plaque containing $\h{z}$ we can take any of the
maps $\pi_{-n}:\h{U}\rightarrow U_{-n}$ with $n>N$. Plaques glue
together to form a Riemann surface $L$. Any such set $L$ will be
called a \textit{leaf} of the regular part.

Once the leaves are endowed with a conformal structure, the map
$\h{f}_c$ becomes a conformal bi-holomorphism between leaves
sending $L(\h{z})$ to $L(\h{f}_c(\h{z}))$.

A result in \cite{LM} states that all leaves in $\mc{R}_c$ are
simply connected. Moreover, if $f$ is a rational function, all the
leaves in $\mc{R}_f$ are simply connected except in the case of
invariant lifts of Herman rings, which are doubly connected. Also,
there are no compact leaves in $\mc{R}_c$. Therefore, by the
Uniformization Theorem, from a conformal view point, leaves are
either disks or planes. The union $\mc{A}_{c}$ of the leaves in
$\mc{R}_c$ conformally isomorphic to the complex plane $\C$ is
called the \textit{affine part} of $\mc{R}_c$. We call any such
leaf an \textit{affine leaf}.

A leaf conformally isomorphic to the disk is called a
\textit{hyperbolic leaf}. Such leaves can exist, for instance,
invariant lifts of Siegel disks belong to hyperbolic leaves.
Examples of non-rotational hyperbolic leaves have been constructed
by Jeremy Khan and Juan Rivera-Lettelier. We can say that the
geometry of the leaves is related to the dynamical properties of
the postcritical set. However, so far there is no characterization
of the existence of hyperbolic leaves (see section 4.1 in
\cite{LM}).

In any case, since the affine part contains all the leaves
containing repelling periodic points, there are always infinitely
many affine leaves in the affine part of any rational map (see
Proposition 4.5 in \cite{LM}). Moreover, for any such leaf $L$
containing a periodic point $\h{p}$, the lift of K\"{o}nig's
coordinate around $\pi(\h{p})$ provides a conformal isomorphism
from $L$ to $\C$.

The following is a simple criterion to check if a leaf $L$ is
affine (see Corollary 4.2 in \cite{LM}). Given a point $\h{z}\in
\mc{R}_c$, if the coordinates $\{z_{-n}\}$ of $\h{z}$ do not converge
to the postcritical set $P(f_c)$ then the leaf $L(\h{z})$
containing $\h{z}$ is affine.

When $c$ is hyperbolic, there are no critical points in the Julia
set, and every backward orbit accumulates on the Julia set, by the
previous argument; if $c$ is hyperbolic then all leaves in
$\mc{R}_c$ are affine.

In order to be an affine lamination, the topology on the affine
part must be such that locally the affine structure on leaves
varies continuously along a transversal. In general, the natural
topology on $\mc{R}_c$ is not enough, since local degree of
plaques may not depend continuously on the transversal direction.

However, for hyperbolic parameters we have:

\begin{proposition}[Lyubich-Minsky]
If $c$ is hyperbolic, the associated regular part $\mc{R}_{c}$ is
a lamination with the topology induced by the natural topology.
\end{proposition}

In Lyubich-Minsky \cite{LM} the concept of convex cocompactness is
introduced in terms of the compactness of the quotient of the
laminated Julia set under dynamics. This justifies the
terminology, however, we use as a definition a proposition in the
same paper which characterizes convex cocompactness.

\begin{definition} A parameter $c$ is called \textit{convex
cocompact} if the critical point is not recurrent and does not
converge to a parabolic cycle.
\end{definition}

In particular, attracting and superattracting parameters are
convex cocompact. If $c$ is convex cocompact, leaves in the
regular part $\mc{R}_{c}$ are all affine, so
$\mc{A}_{c}=\mc{R}_{c}$; also, the affine part in this case has a
simple description: $\mc{R}_{c}=\mc{N}_{c} \setminus
\{\textnormal{attracting and parabolic cycles}\}$ (see Proposition
4.5 in \cite{LM}).

\section{Topology of inverse limits}

\subsection{The laminated Julia set}

Let $\mc{J}_c=\pi^{-1}(J_c)\cap \mc{R}_c$ be the lift of the Julia
set on the regular part. We call $\mc{J}_c$ the \textit{laminated
Julia set} associated to $f_c$. If $c$ is a convex co-compact
parameter, by a theorem in \cite{LM}, $\mc{J}_c$ is compact inside
$\mc{R}_c$ with the natural topology, so in this case
$\mc{J}_c=\pi^{-1}(J_c)$.

In this section we discuss the topological properties of
$\mc{J}_c$. Given a repelling periodic point $q$ in the dynamical
plane, let us denote by $\h{q}$ the periodic point in the regular
part which satisfies $\pi(\h{q})=q$. Let $\mathcal{P}$ be the set
of repelling periodic points in $\mc{R}_{c}$.

\subsubsection{Leafwise connectivity of $\mc{J}_c$}

\begin{lemma}\label{lcdisconnection}
Let $c$ be a parameter with locally connected Julia set,
$\overline{P(f_c)}\neq J_c$, and such that $\pi^{-1}(J(f_c))$
contains an irregular point in the natural extension $\mc{N}_c$.
Then, there is a leaf $L$ in the regular part $\mc{R}_{c}$, such
that $\mc{J}_c \cap L$ is disconnected.
\end{lemma}

\begin{proof} Let $\hat{z}$ be an irregular point in the Julia set in
$\mc{N}_c$ then, since $J(f_c)$ is completely invariant and
locally connected, $\pi_{-j}(\hat{z})=z_{-j}$ belongs to $J(f_c)$
and there is an external ray $R(j)$ landing at $z_{-j}$ for every
$j$. Let $r_0$ be a point in $R(0) \cap A(\infty)$, by pulling
back $r_0$ along the backward orbit determined by $\hat{z}$, there
is a point $\h{r} \in \mc{R}_c$ such that $\pi(\h{r})=r_0$ and
$\pi_{-j}(\h{r})\subset R(j)$, now by moving $r_0$ along the ray
$R(0)$, we construct a line $\h{R}$ in the regular part, such that
$\pi_{-j}(\h{R})=R(j)$. Let $L$ be the leaf in the regular part
containing $\h{R}$. By construction the endpoints of $\h{R}$ are
the irregular points $\h{z}$ and $\h{\infty}$, this line can not
have accumulation points in $\mc{R}_c$ when $r_0$ either tends to
$z_0$ or to $\infty$. So, $\h{R}$ is a line escaping to infinity
in both directions and separates $L$ in two pieces.

Now, fix $r_0$ in the external ray $R(0)$, let $a\in J(f_c)
\setminus\overline{P(f_c)}$, and choose two paths, $\sigma_1$ and
$\sigma_2$, from $r_0$ to $a$ starting at different directions
with respect to the ray $R(0)$, and such that none of them crosses
$R(0)$ again. These two paths lifts to paths in $L$, joining
$\h{r}$ with points in $\mc{J}_c\cap L$, by construction the
points lie on different sides of the line $\h{R}$. Therefore
$\mc{J}_c\cap L$ is disconnected.
\end{proof}

As interesting examples of polynomials with $\mc{J}$ non leaf-wise
connected, consider: quadratic polynomials with parabolic cycles
(see Tomoki Kawahira's paper \cite{Tom}), or the Feigenbaum
quadratic polynomial. In the case of quadratic polynomials with
parabolic cycles, leaves in the regular part where the Julia set
is disconnected are precisely the periodic leaves corresponding to
the parabolic cycle. As for repelling cycles the linearizing
coordinate, in this case Fatou's coordinate, gives a
uniformization of the periodic leaves on the regular part. A
detailed explanation and construction of such uniformization can
be found in \cite{Tom}. It is worth noting that, when $f_c$ is
parabolic, then for any leaf $L$ the set $\mc{J}_c\cap L$ consists
of finitely many components.

In general, it is not clear whether the intersection of the
laminated Julia set consists of finitely many pieces. An
interesting example to bear in mind is the natural extension
associated to the Feigenbaum polynomial. In this case, there are
irregular points in the fiber of the Julia set, and then there are
leaves for which the Julia set is not path connected.

As a corollary, if under the conditions of Lemma
\ref{lcdisconnection} the laminated Julia set $\mc{J}_c$ is
leaf-wise connected, then there are no irregular points in
$\mc{J}_c$, which implies that $f_c$ is convex co-compact. The
following Lemma shows that the converse is also true:

\begin{lemma}\label{leafwise.con} If $c$ is a convex cocompact
parameter, then the set $\mc{J}_c$ is leafwise connected.
\end{lemma}

\begin{proof} Since $f_c$ is convex cocompact, then the Julia set
$J(f_c)$ is locally connected and every external ray lands in the
Julia set. Moreover, since there are no irregular points in
$\pi^{-1}(J(f_c)$ in $\mc{N}_{f_c}$, every solenoidal ray also
lands at $\mc{J}_c$. Let $L$ be a leaf in the regular part, and
let $\h{z}$ and $\h{\zeta}$ be two points in $L\cap \mc{J}_c$.
Since $\h{z}$ and $\h{\zeta}$ belong to the same leaf, there exist
a path $\sigma$ in $L$ joining $\h{z}$ with $\h{\zeta}$. Then by
using the external ray flow in $L$ the path $\sigma$ projects onto
a path whose trajectory lies in $\mc{J}_c$ and connects $\h{z}$
with $\h{\zeta}$, so $\h{z}$ and $\h{\zeta}$ lie in the same
connected component of $L\cap \mc{J}_c$.
\end{proof}

Together, Lemma~\ref{lcdisconnection} and Lemma~\ref{leafwise.con}
imply:

\begin{proposition}
A quadratic parameter $c$ with locally connected Julia set
$J(f_c)$ and $\overline{P(f_c)}\neq J_c$ has leafwise connected
Julia set $\mc{J}_c$ in its regular part $\mc{R}_c$, if and only
if $c$ is convex cocompact.
\end{proposition}

\subsection{Hyperbolic components}

As we pointed on before, combinatorially, the Julia sets
associated to parameters within a hyperbolic component $H$ are
indistinguishable from that associated to the root of $H$, . Here,
we present a proposition that topologically ties the regular parts
of parameters inside hyperbolic components with the regular part
of their center.

Although it was not explicitly stated the proof of the next
proposition follows immediately from Lemma 11.1 in Lyubich and
Minsky. We include the settings and the statement of that lemma
for reference. The interested reader can find the proof in
\cite{LM}.

Let $U$ and $V$ be two open sets with $\bar{U}\subset V$ and let
$f:U\rightarrow V$ be an analytic branched covering map. Let us
remark that when $U$ and $V$ are disks with the property $U\subset
\subset V$ the map $f$ is called \textit{polynomial-like} or
\textit{quadratic-like} if the degree of $f$ is $2$. Let
$\mc{N}_f$ denote the set of backward orbits of $f$. In this
setting, iterations of the map $\h{f}$ on $\mc{N}_f$ may not be
defined, since $f$ is not defined on $V\setminus U$. For
$m=1,2,...$, let $\mc{N}_m\subset \mc{N}_f$ be the set of backward
orbits $\h{z}$ that can be iterated under $\h{f}$ \textit{at most}
$m$ times. So, we have inclusions $\mc{N}_m\subset \mc{N}_{m+1}$
for $m=1,2,...$.

The map $\h{f}^{-m}:\mc{N}_f \rightarrow \mc{N}_f$ is an
immersion, that maps $\mc{N}_f$ onto $\mc{N}_{m}$ for $m=1,2,...$.
So, by composing with the inclusions $\mc{N}_{m}\hookrightarrow
\mc{N}_f$, we consider $\mc{N}_f$ as an extension of $\mc{N}_{m}$
and denote these extensions by $\mc{N}_m$. Make
$\mc{N}_f=\mc{N}_0$ and identify any point $\h{z}\in \mc{N}_m$
with $\h{f}^{-1}(\h{z})$, so the map $\h{f}^{-1}$ induces the
following increasing sequence of sets
$$\mc{N}^0\hookrightarrow\mc{N}^1\hookrightarrow \mc{N}^2\hookrightarrow...$$
let $\mc{D}_f=\cup \mc{N}_m$, a set $W$ in $\mc{D}_f$ is said to
be open if $W\cap \mc{N}_m$ is open for every $m$. The set
$\mc{D}_f$ is called the \textit{direct limit} of the increasing
sequence above. The natural extension $\h{f}$ of $f$ in $\mc{N}_f$
induces a homeomorphism of $\mc{D}_f$ into itself. Now, we can
state the following Lemma:

\begin{lemma}[Lyubich and Minsky]\label{lem.11.1}
Assume that a branched covering $f:U\rightarrow V$ is the
restriction of a rational endomorphism $R:\bar{\C}\rightarrow
\bar{\C}$ such that $\C\setminus V$ is contained in the basin of
attraction of a finite attracting set $A$. Then
$\h{f}:\mc{D}_f\rightarrow \mc{D}_f$ is naturally conjugate to
$\h{R}:\mc{N}_R\setminus \h{A} \rightarrow \mc{N}_R \setminus
\h{A}$.
\end{lemma}

Let $c$ be a hyperbolic parameter contained in the hyperbolic
component $H_c$ in the Mandelbrot set. Then we have the following:

\begin{proposition}\label{attract.sup} The regular part
$\mc{R}_c$ is homeomorphic to the regular part of the quadratic
polynomial associated to the center of $H_c$.
\end{proposition}

\begin{proof} First, let us discuss parameters inside the Main Cardioid.

Let $f_{\epsilon}(z)=z^2+\epsilon$ a quadratic polynomial with
$\epsilon$ a parameter inside the Main Cardioid. Thus
$f_{\epsilon}$ has an attracting fixed point $a_{\epsilon}$. The
Fatou set consists of two open sets corresponding to the basins of
infinity and $a_\epsilon$, $A(\infty)$ and $A(a_\epsilon)$. The
Julia set $J(f_{\epsilon})$ is a quasicircle, so that the
conjugating map $\phi:J(f_\epsilon)\rightarrow \mathbb{S}^1$ can
be quasiconformally extended to a neighborhood of $J(f_\epsilon)$.
Since $f_\epsilon$ is expanding on the Julia set, the map
$\phi\circ f_\epsilon \circ \phi^{-1}$ is an expanding circle map
of degree 2.

By a Theorem of Shub the map $\phi\circ f_\epsilon \circ
\phi^{-1}$ is topologically conjugate to $z^2$, that is, there is
a map $h:\mathbb{S}^1\rightarrow \mathbb{S}^1$ such that
$f_0=h\circ \phi\circ f_\epsilon \circ \phi^{-1} \circ h^{-1}$.
Also, $h$ admits an equivariant extension to a neighborhood of
$\mathbb{S}^1$. Actually, B\"ottcher's coordinate in the basin of
infinity extends $h$ to the whole basin of infinity. So we obtain
a conjugacy of $f_\epsilon$ to $f_0$ defined on a simply connected
neighborhood $U$ containing the basin of infinity and the Julia
set $J(f_\epsilon)$. We can choose $U$ small enough such that the
map $f_\epsilon :U\ra V$ is quadratic-like. This implies that the
map $f_\epsilon:U\rightarrow V$ is topologically equivalent to the
map $f_0:\phi(U)\rightarrow \phi(V)$. By construction,
$\C\setminus V$ is contained in the basin of attraction of
$a_\epsilon$.

By Lemma~\ref{lem.11.1}, the map
$\h{f}_\epsilon:\mc{D}_{f_\epsilon}\rightarrow\mc{D}_{f_\epsilon}$
is conjugate to $\h{f}_\epsilon:\mc{N}_{\epsilon}\setminus
\h{a}_\epsilon \rightarrow \mc{N}_{\epsilon}\setminus
\h{a}_\epsilon$, also the map $\h{f}_0:\mc{D}_{f_0}\rightarrow
\mc{D}_{f_0}$ is conjugate to $\h{f}_0:\mc{N}_0\setminus \h{0}
\rightarrow \mc{N}_{0}\setminus \h{0}$. But the conjugacy $\phi$
from $f_\epsilon$ to $f_0$ lifts to a conjugacy from
$\h{f}_\epsilon:\mc{D}_{f_\epsilon}\rightarrow\mc{D}_{f_\epsilon}$
to $\h{f}_0:\mc{D}_{f_0}\rightarrow\mc{D}_{f_0}$. So,
$\mc{R}_\epsilon=\mc{N}_\epsilon \setminus
\{\h{\infty},\h{a}_\epsilon\}$ is homeomorphic to
$\mc{R}_0=\mc{N}_0 \setminus \{\h{\infty},\h{0}\}$.

It follows that $\mc{N}_\epsilon$ is homeomorphic to the double
cone over $\sol$. Let us remark that since the homeomorphism above
is a conjugacy, it sends the Julia set $\mc{J}_\epsilon$ onto
$\sol=\mc{J}_0$, so we can restrict such homeomorphism to the lift
of the basin of attraction of $a_\epsilon$, therefore
$\inver(f_\epsilon, A(a_\epsilon))\cup J(f_\epsilon)$ is
homeomorphic to $\cone$.

Now, let $f_c(z)=z^2+c$ be a quadratic polynomial with $c$ in a
hyperbolic component $H$, with attracting cycle
$P=\{p_1,...,p_n\}$ of period $n$. Let $U_1,U_2,...,U_n$ be the
Fatou components of the basin of attraction of $P$.

The regular part $\mc{R}_c$ can be decomposed in several parts by
cutting along the Julia set $\mc{J}_c$, which are:

\begin{itemize}

\item The Julia set $\mc{J}_c$.

\item The lift of the basin at infinity $\pi^{-1}(A(\infty))$ in
$\mc{R}_c$.

This set is homeomorphic to $\sol \times (0,1)$ by means of the
lift of B\"ottcher's coordinate.

\item Fatou components of finite branching.

This set is homeomorphic to a countable union of sets
$\{\mc{V}_n\}$, where each $\mc{V}_n$ is homeomorphic to
$\mathbb{D}\times \cantor$. Each point $\h{z}\in \mc{V}_n$ has at
most finitely many coordinates in the immediate basin of
attraction of the critical cycle $\cup U_j$. Also, each component
in $\mc{V}_n$ projects onto some $U_i$ for $i$ fixed.

\item The invariant lift $\h{U}$ of the basin of attraction of
$P$.

Let $\h{U_i}$ be the of points $\h{z}\in\h{U}$ such that
$\pi(\h{z})\in U_i$, with the lift of K\"onig's coordinate on
$\h{U}_i$, $\h{f}^n_c$ has the form of $z\mapsto \lambda z$, where
$\lambda$ is the multiplier of the cycle $P$, so $\h{f}^n_c$ in
$\h{U}_i$ is topologically equivalent to $\h{f}_\lambda$ in
$\pi^{-1}(A(a_{\lambda}))$. By the discussion above, $\h{U}_i$ is
homeomorphic to $\sol\times (0,1)$.

\end{itemize}

Let $c_0=H(0)$ be the center of $H$. Now, the regular part
$\mc{R}_{c_0}$ has the same decomposition as $\mc{R}_c$. This
decomposition is such that the corresponding components are
homeomorphic. These homeomorphisms glue together to a
homeomorphism from $\mc{R}_c$ to $\mc{R}_{c_0}$. \end{proof}

We call this decomposition of $\mc{R}_c$, the \textit{laminated
decomposition} of the regular part associated to the hyperbolic
parameter $c$. Tomoki Kawahira independently proved
Proposition~\ref{attract.sup} in the more general setting of
hyperbolic rational maps and in the quasiconformal level. That is,
hyperbolic affine laminations are qc-stable in the Lyubich-Minsky
setting (see \cite{Tom2}). Let $H$ be a hyperbolic component, and
$c_0=H(0)$ be the center of $H$. For any hyperbolic parameter $c$
in the boundary of $H$, the path from $c_0$ induces a
transformation $h_c$ from $\mc{R}_{c_0}$ to $\mc{R}_c$. If
correspondingly $c_1$ denotes the root of $H$, then $f_{c_1}$ is a
parabolic quadratic polynomial. Let $\mc{L}$ be the regular part
of $\mc{R}_{c_0}$ with the leaves containing the dynamic root
cycle removed. Analogously, let $\mc{L}'$ be set obtained by
removing from $\mc{R}_{c_1}$ the periodic leaves associated to
Fatou's coordinate. Then another result of Kawahira (see
\cite{Tom}), states that the map $h_{c_1}=\lim h_c$, as $c$ tends
to $c_1$ in $H$, is a laminar homeomorphism between $\mc{L}$ and
$\mc{L}'$, and moreover, the map $h_{c_1}$ semi-conjugates
$\h{f}_{c_0}|_{\mc{L}}$ to $\h{f}_{c_1}|_{\mc{L}'}$.

\medskip

Let $\mc{F}_c=\pi^{-1}(F(f_c))\cap\mc{R}_{c}$ be the lift of the
Fatou set to the regular part; we call $\mc{F}_c$ the
\textit{laminated Fatou set}. Given an affine leaf $L$ in
$\mc{R}_{c}$, consider the uniformization $\phi:L\rightarrow \C$.
We call a subset $A$ in $L$ \textit{bounded} if the corresponding
set $\phi(A)$ is bounded in $\mathbb{C}$. As a consequence of the
laminated decomposition of regular parts in the proof of the
previous lemma we have the following:

\begin{lemma} Let $f_c$ be a convex cocompact quadratic polynomial,
and let $A$ be a connected component of $\mc{F}_c$ inside an
affine leaf $L$. Then $A$ is bounded if and only if the
restriction $\pi|_A$ has finite degree.
\end{lemma}

\begin{proof}

If $c$ is not hyperbolic, then in the dynamical plane the Fatou
set $F(f_c)$ consist only of the basin of infinity. The lift of
$F(f_c)$ consist just of the solenoidal cone at infinity; clearly
all Fatou components $A$ of $\mc{F}_c$ are unbounded and the
degree of $\pi$ restricted to $A$ is infinity.

When $c$ is hyperbolic, then by Lemma~\ref{attract.sup} we can
assume that $c$ is superattracting. In this case, all unbounded
Fatou components belong either to the solenoidal cone at infinity
or to the solenoidal cones at the critical orbit. The degree of
$\pi$ restricted to each of these Fatou components is infinity.
The remaining Fatou components, the bounded ones, belong by
definition to components where the degree of $\pi$ is finite.
\end{proof}

\subsection{Unbounded Fatou components}

We want to count how many unbounded Fatou components there are in
$L$. Recall that the valence $v_p$ of a repelling periodic point
$p$ is the number of external rays landing at $p$. When $p$ is the
dynamic root point $r_c$ and $c$ belongs to a satellite hyperbolic
component, there are $v_p$ Fatou components touching at $r_c$. If
$c$ belongs to a primitive hyperbolic component, $r_c$ only
touches one Fatou component and $v_{r_c}=2$, see Milnor's
\cite{Mext}. As noted in Lyubich and Minsky \cite{LM}, the
K\"{o}nigs's coordinate $\phi$ around a repelling periodic point
$p$ lifts to the uniformization $\Phi$ of $L(\h{p})$. Moreover, if
$m$ is the period of $p$, then $\Phi$ conjugates
$\h{f}_c^m|_{L(\h{p})}$ to the affine map $z\mapsto \lambda_p z$,
where $\lambda_p$ is the multiplier of $p$. By means of $\Phi$,
local properties of $p$ are reflected in global properties in
$L(\h{p})$. This is the idea behind the proof of the following
proposition.

\begin{proposition}\label{unbounded} Let $f_c$ be a quadratic
polynomial and let $L$ be a periodic affine leaf in the regular
part $\mc{R}_{f_c}$ containing a periodic point $\h{p}$. Then, the
number of unbounded Fatou components of $L\setminus \mc{J}_c$ is
either $v_p$ if the corresponding periodic point $p$ does not
belongs to the dynamic root cycle, and otherwise there are two
cases: There are $2v_p$ if $c$ belongs to a satellite component
and exactly 3 unbounded Fatou components if $c$ belongs to a
primitive hyperbolic component.
\end{proposition}

\begin{figure}[htbp]
\begin{center}
\input{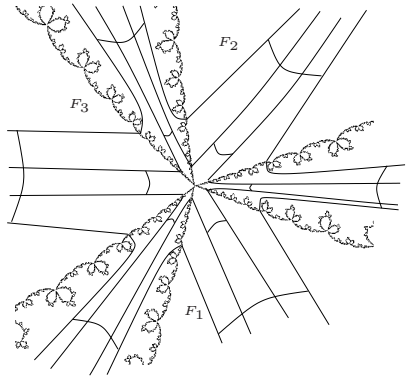}
\caption{The annulus $A$ when $p$ belongs to the dynamic root
cycle.}
 \label{unbounded2}
\end{center}
\end{figure}

\begin{proof} Let $m$ be the period of $\h{p}$; then $p=\pi(\h{p})$ also has
period $m$. Since Siegel periodic points lift into hyperbolic
leaves in $\mc{R}_c$, by Lemma~\ref{lamin.julia}, $p$ most be a
repelling periodic point and so, it belongs to the Julia set
$J(f_c)$. In the dynamical plane, there are $v_p$ rays landing at
$p$ that cut $\C$ in $v_p$ sectors. Let $S$ be one of these
sectors, by construction $S$ is invariant under an appropriate
iterate of $f_c$, say $f_c^k$, where $k$ is a multiple of $m$.

Every ray landing at $p$ lifts to a landing ray at $\h{p}$. As in
the dynamical plane, landing rays cut the leaf $L(\h{p})$ in $v_p$
sectors. We will check first that when $p$ is not in the dynamic
root cycle, every ray landing at $\h{p}$ determines an unbounded
component of the Fatou set in $L(\h{p})$, to do so we prove that
every sector $\h{S}$ in $L(\h{p})$ intersects the Fatou set in two
unbounded components. Since $\h{f}_c^k$ is a similarity in
$L(\h{p})$, it suffices to prove that the intersection of the
Julia set $\mc{J_c}$ with a sector $\h{S}$ grows to infinity in
one direction.

To do that, let us construct a fundamental piece to the action of
$f_c^k$ in $\C$. Given a sector $S$ in $\C$, let $E_S=J(f_c)\cap
S$. Let us take a point $b$ in $E_S$, very close to $p$, and such
that there is a pair of rays $R_{b}$ and $R'_{b}$ landing at $b$,
whose images $R_{f^k_c(b)}$ and $R'_{f^k_c(b)}$ belong to the wake
$W_S$ determined by $R_b$ and $R'_b$. Fix an equipotential $E_r$
and join consecutive landing rays by arcs of this equipotential.
We obtain a region $P$ around $p$, since $p$ is repelling the set
$P$ is compactly contained in $f_c^k(P)$ and the annulus
$A=\overline{f_c^{k}(P)\setminus P}$ is the fundamental piece we
are looking for (see Figure~\ref{fig.unbounded1}).

\begin{figure}[htbp]
\begin{center}
\input{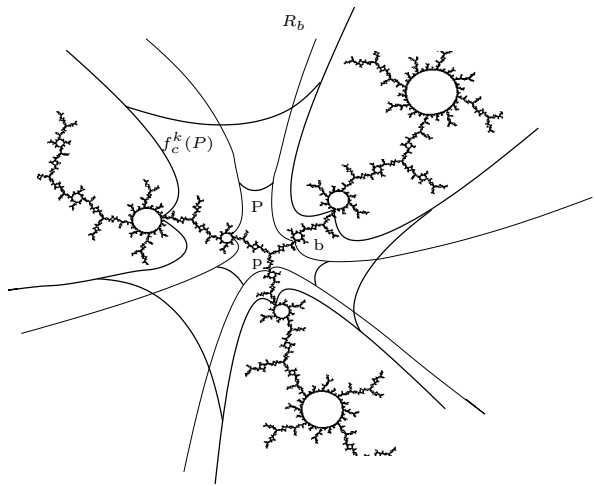}
\caption{The annulus $A$ when $p$ is not in the dynamic root
cycle.} \label{fig.unbounded1}
\end{center}
\end{figure}

The Julia set intersects the annulus $A$ at the wake $W_S$ defined
by the rays $R_b$ and $R'_b$, thus for every $S$, we can enclose
the subset of the Julia set in $A\cap W_S$, in a simply connected
open set $V_S$ contained in $A\cap W_S$.

Now, $A$ lifts to an annulus $\h{A}$ in $L$, which by
construction, is a fundamental region for the action of
$\h{f}_c^k$. By iterating $\h{f}_c^k$ on $\h{V}_S$, the lifts of
the sets $V_S$, we obtain a set in $L$, similar to a ``string of
pearls" (see Figure~\ref{fig.unbounded1b}), which contains the
subset of the Julia set $\h{E}$ in $\mc{J}_c$. By construction,
the lift of this ``string of pearls" into $L(\h{p})$ is unbounded
and there is one for each sector $S$ in the dynamical plane. The
conclusion follows.

\begin{figure}[htbp]
\begin{center}
\input{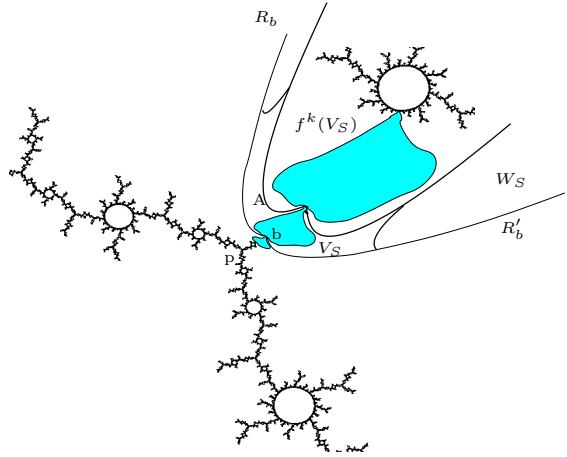}
\caption{The set $V_S$ and some iterates.} \label{fig.unbounded1b}
\end{center}
\end{figure}

In the case when $p$ belongs to the the dynamic root cycle, there
are extra unbounded components of the laminated Fatou set in
$L(\h{p})$ coming from the lifts of the Fatou components attached
to $p$. We have to modify the previous argument to sectors that
contain a Fatou component attached to $p$. Let $S$ be such a
sector containing a Fatou component, say $F_S$, having $p$ in the
boundary. Instead of the point $b$ as above, we consider a pair of
points $q$ and $q'$ on $\partial F_S$ and on opposite sides of
$p$. There are two pairs of landing rays, $\{R_q,I_q\}$ and
$\{R_{q'},I_{q'}\}$ landing at $q$ and $q'$. The $R$'s are
external and the $I$'s are internal rays. On the basin of infinity
we follow the same construction as in the case before, whereas in
$F_S$ connect the internal rays $I_q$ with $I'_q$ by a internal
equipotential. So we obtain again a puzzle piece $P$ around $p$,
see Figure \ref{unbounded2}. From now on, the argument above goes
through either for $S$ or sectors that do not contain Fatou
components attached to $p$. To complete the computation, notice
that if $c$ belongs to a satellite component, each pair of
consecutive rays landing at $p$ contains a Fatou component with
$p$ on its boundary, so there are $2v_p$ components. In the
primitive case, there are only 2 rays landing at $p$ and only 1
Fatou component contains $p$ on its boundary.
\end{proof}

The following proposition gives another restriction that
combinatorics impose on the leaf structure of the Julia set.

 \begin{proposition} Let $f_c$ be a convex cocompact
quadratic polynomial, and let $L$ be a non-periodic affine leaf.
Then, the number of unbounded Fatou components of $L\setminus
\mc{J}_c$ is either $1$ or $2$.
\end{proposition}

\begin{proof} Let $R$ and $R'$ be two external rays landing at the
Julia set such that the wake $W$ determined by $R$ and $R'$ does
not contain postcritical points. If $\h{R}$ and $\h{R}'$ are lifts
of the rays $R$ and $\h{R}$ landing at the same point in a leaf
$L$ in $\mc{R}_c$, then every arc connecting $R$ and $R'$ along
$W$ must lift to an arc joining $\h{R}$ and $\h{R}'$. Such an arc,
for instance, can be taken to be part of an equipotential. Thus,
if $\h{W}$ is the wake determined by $\h{R}$ and $\h{R}'$, then
$\h{W}\cap \mc{J}_c$ is a bounded set in $L$. So, $\h{R}$ and
$\h{R}'$ belong to the same unbounded Fatou component.

Now, let $L$ be a non-periodic affine leaf. Since the Julia set
$J(f_c)$ is locally connected, the set $\mc{J}_c$ is leafwise
locally connected and every point in $\mc{J}_c$ is the landing
point of some external ray in $\mc{R}_c$. If $L$ has more than $3$
unbounded Fatou components, by Lemma~\ref{leafwise.con},
$\mc{J}_c\cap L$ is path connected, so there is a point $\h{z}$ in
$\mathcal{J}_c \cap L$ on the boundary of at least 3 unbounded
Fatou components in $L$. This implies that $\h{z}$ is the landing
point of at least $3$ rays, each ray in a different unbounded
Fatou component. Hence, each coordinate $\pi_{-n}(\h{z})\in
J(f_c)$ is also the landing point of at least 3 rays. By the
argument above, the rays landing at $\pi_{-n}(\h{z})$ cut the
postcritical set in three disjoint pieces. This implies that
$\pi_{-n}(\h{z})$ must be a vertex of the Hubbard tree of $J(f_c)$
for each $n$. Since $L$ is non-periodic, $\h{z}$ is non-periodic,
and the set of coordinates $\pi_{-n}(\h{z})$ is an infinite set in
the dynamical plane. But, this contradicts the fact that the set
of vertices in a Hubbard tree with degree at least $3$ is finite.
\end{proof}

Let us note that when $L$ is a non-periodic affine leaf, the cases
where $L$ has $1$ or $2$ unbounded Fatou components can both
happen. To get leaves with two unbounded Fatou components,
construct a backward history of \textit{biaccesible} points
$\h{z}$ in $\mc{J}_c$, that is, points where at least 2 rays land,
with the property that the two wakes determined by the rays
landing at $z_{-n}$ contain postcritical points. The case where
$L$ has only $1$ unbounded component is the most common inside
regular parts of convex cocompact parameters $c$, because almost
all repelling periodic points in $J(f_c)$ are the landing point of
exactly one ray. For instance, the leaf containing the lift of the
$\beta$ fixed point is always a leaf with exactly one unbounded
component. Non-periodic leaves with one unbounded Fatou components
can also be constructed.

In regular parts of convex cocompact parameters, leaves with more
than three unbounded components are in one-to-one correspondence
with the vertexes of the Hubbard tree with degree greater than 2.
All other leaves either have one or two unbounded components. So,
there are only finitely many leaves with more than three unbounded
components. This is how the number of unbounded Fatou components
is related to the combinatorics of the parameter $c$. The
following proposition describes the cycle of leaves with more than
3 unbounded components using internal addresses.

\begin{proposition}\label{internal.addr} Let $c$ be
a superattracting parameter with internal address $1\rightarrow
n_1 \rightarrow n_2 \rightarrow ... \rightarrow n_k$. If
$n_{j-1}\mid n_j$,  for $j<k$ then there is a cycle of $ n_{j}$
periodic leaves such that each leaf has $\frac{n_j}{n_{j-1}}$
unbounded Fatou components.  When $j=k$ there are $n_{k}$ leaves
with $2 \frac{n_k}{n_{k-1} }$ unbounded Fatou components. If
$n_{j-1}\nmid n_j$, there is a cycle of $n_{j}$ periodic leaves
with 2 unbounded Fatou components for $j<k$, and 3 unbounded Fatou
components if $j=k$. \end{proposition}

\begin{proof} The condition of whether $n_{j-1}$ divides $n_j$ or not
reflects whether the combinatorial arc of $c$ crosses a satellite
or a primitive hyperbolic component in the parameter plane. Let
$j=k$, the number $n_k$ corresponds to the period of the critical
orbit, which is equal to the period of the dynamic root point;
hence by Proposition~\ref{unbounded}, if $n_{k-1}\mid n_k$ the
parameter belongs to a satellite component, in which case the
valence of the dynamic root point is $\frac{n_k}{n_{k-1}}$, and
there are $2 \frac{n_k}{n_{k-1}}$ unbounded Fatou components. If
$n_{k-1}\nmid n_k$, then $c$ is the center of a primitive
hyperbolic component and there are $n_k$ leaves with 3 unbounded
Fatou components.

Now, let $j<k$. If $n_{j-1}\mid n_j$, then there is a repelling
cycle $P$ of $f_c$ with period $n_j$ and valence
$\frac{n_j}{n_{j-1}}$, and by Proposition~\ref{unbounded} the lift
$\h{P}$ belongs to a cycle of $n_j$ periodic affine leaves with
$\frac{n_j}{n_{j-1}}$ unbounded Fatou components in
$\mc{R}_{f_c}$. If $n_{j-1}\nmid n_j$, the corresponding ray
portrait has valence $2$ and period $n_j$.
\end{proof}

\begin{figure}[htbp]
\begin{center}
\input{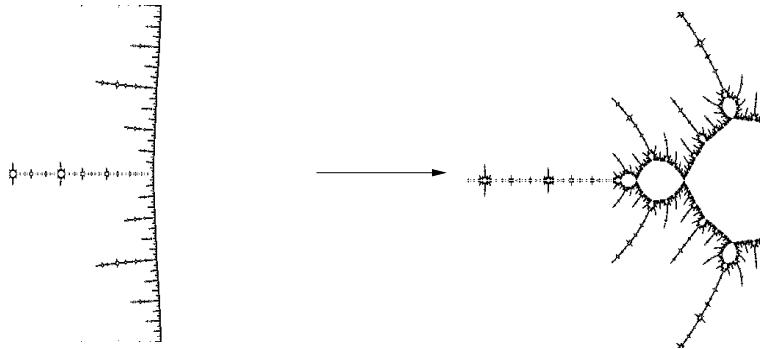}
\caption{Collapsing of unbounded Fatou components by bifurcation.}
 \label{airbif}
\end{center}
\end{figure}

Leaves containing the periodic lift of the dynamic root cycle of
primitive parameters have exactly 3 unbounded Fatou components. If
we change the parameter $c$ to any center in a adjacent satellite
component, one of the unbounded Fatou components collapses to an
infinite number of \emph{bounded} Fatou components. Although $c$
still has a cycle corresponding to the dynamic root cycle of the
previous satellite parameter, the leaves containing this cycle
will have only $2$ unbounded Fatou components, each of them on the
lift of the basin of infinity $\pi^{-1}A_c(\infty)$ (see
Figure~\ref{airbif}).

With all the previous discussion, now we can prove a special case
of the Main Theorem, namely when the homeomorphism $h$ is a
conjugacy:

\begin{proposition}
Let $c_1$ and $c_2$ be two superattracting parameters. If
$h:\mc{R}_{c_1} \rightarrow \mc{R}_{c_2}$ is a homeomorphism
conjugating $\h{f}_{c_1}$ with $\h{f}_{c_2}$, then $c_1=c_2$.
\end{proposition}

\begin{proof} Any such conjugation sends periodic points into
periodic points. Hence, by Lemma~\ref{lamin.julia}, $h$ sends the
Julia set $\mc{J}_{c_1}$ into the Julia set $\mc{J}_{c_2}$. Since
$h$ is a homeomorphism it has to leave invariant the number of
unbounded Fatou components, and as $h$ is a conjugacy of dynamics,
it also leave invariant the combinatorial rotation numbers among
the unbounded Fatou components on periodic leaves. This means that
$c_1$ and $c_2$ must have the same labelled internal address. But
by Lau and Schleicher this implies $c_1=c_2$.
\end{proof}



Given a repelling periodic point $q$ in the dynamical plane, let
us denote by $\h{q}$ the periodic point in the regular part which
satisfies $\pi(\h{q})=q$. Let $\mathcal{P}$ be the set of
repelling periodic points in $\mc{R}_{c}$.

\begin{lemma}\label{lamin.julia}
The set $\mathcal{J}_c$ is a closed and perfect set in
$\mc{R}_{c}$. Every periodic point in $\mc{R}_{c}$ either belongs
to $\mc{P}$ or is a Siegel periodic point. Moreover,
$\mc{J}_c=\bar{\mc{P}}$.
\end{lemma}

\begin{proof}  The projection $\pi$ is a continuous function from
$\mc{R}_{c}$ to $\mathbb{C}$. Since $J(f_c)$ is closed and
perfect, the set $\mc{J}_c$ inherits these properties from the
dynamical plane.

The fact that lifts of Siegel periodic points in the dynamical
plane belong to the regular part is a consequence of the existence
of linearizing coordinates around Siegel periodic points. By a
Theorem of Fatou, parabolic and attracting cycles are the
accumulation set of some critical orbit, hence the pull back
$\{U_{-n}\}$ of every neighborhood $U_0$ contains critical points
at infinitely many times $n$, so parabolic and attracting cycles
lift to irregular points.

Now, let us prove that Cremer cycles also lift to irregular
points. This is, actually, a consequence of the Shrinking Lemma.
To illustrate its use, we will follow the proof given by Lyubich
and Minsky \cite{LM}, see also Proposition 1.10 in Lyubich's
survey \cite{L}.

Suppose on the contrary, that a Cremer cycle lifts to a periodic
regular point. By considering an appropriate iterate of $f_c$, we
can assume that the cycle is a Cremer fixed point $a_0$. As
$\h{a}_0$ is regular, there exist an open neighborhood $U_0$ of
$a_0$ such that no component $U_{-n}$ of the pull back of $U_0$
along $\h{a}_0$ contains critical points. Without loss of
generality we can assume that $U_0$ is a small disk around $a_0$,
then $U_{-n}$ is also a topological disk, since
$f:U_{-n}\rightarrow U_{-n+1}$ is a conformal isomorphism. Hence,
there is a sequence of Riemann maps $\phi_n:\mathbb{D}\rightarrow
U_{-n}$ with $\phi_n(0)=a_0$. Put
$\rho_n=f^n\circ\phi_n:\mathbb{D}\rightarrow U_0$, then by
Montel's Theorem, the sequence $\{\rho_n\}$ is a normal family.
Let $P=\{\alpha_1,\alpha_2,...,\alpha_p\}$ be any cycle of $f_c$
of period $p\geq 3$. By normality there is a disk
$D'\subset\subset \mathbb{D}$ around $0$ such that $\phi_n(D')$
does not intersect $P$ for every $n$. Therefore,
$\{\phi_n|_{D'}\}$ is also a normal family.

Let $B(a_0,\delta)$ denote the ball around $a_0$ of radius
$\delta$, since $\rho_n$ is a conformal isomorphism, for $\delta$
small enough $\textnormal{mod} (\rho^{-1}_n(U_0 \setminus
B(a_0,\delta))=\textnormal{mod}(U_0\setminus B(a_0,\delta))$, so
the modulus of $\rho^{-1}_n(U_0 \setminus B(a_0,\delta))$ depends
only on $\delta$. Because $\textnormal{mod}(U_0\setminus
B(a_0,\delta)) \rightarrow \infty$ when $\delta \rightarrow 0$,
there exist a $\delta$ such that $B=B(a_0,\delta)\subset
\rho_n(D')$ for all $n$.

It follows that $\textnormal{diam}(\phi^{-1}_n(B))\rightarrow 0$
uniformly when $n$ tends to infinity, otherwise by taking a
converging subsequence from $\{\phi_n\}$ and by normality there
would be a limiting open set $B_\infty$ containing $a_0$ such that
$f^{n_k}(B_\infty)\subset U_0$ for all $k$, contradicting the fact
that every Cremer periodic point belongs to the Julia set
$J(f_c)$. If the diameters of $B_n$ tend to zero, there is an $m$
such that $B_m$ is compactly contained in $B$ and $f^m(B_m)=B$,
but this would imply that $f$ is repelling at $a_0$, again a
contradiction.

Take $\hat{z}\in \mathcal{J}_c$, by continuity of $f_c$ and the
density of the set of repelling periodic points in $J(f_c)$, there
is a periodic point $p$ in the $\delta$ neighborhood of $z_{-n}$
such that $|f_c^j(p)-z_{-n+j}|<\epsilon$ for $j=0,...,n$. Then,
$\widehat{f_c^n(p)}$ is a repelling periodic point in
$B(D(z_0,\epsilon),\hat{z},n)$.
\end{proof}

As part of the proof of the previous lemma let us remark the
following:

\begin{corollary} Every point in the invariant lift of a Cremer
cycle is irregular.
\end{corollary}

\subsection{Topology of ends}

In this section, we will consider the regular part $\mc{R}_c$
endowed with the natural topology, so that $\mc{R}_c$ is locally
compact. Topology is important for local compactness since when we
endow affine laminations with Lyubich-Minsky topology, Lyubich and
Lasse Rempe [Lyu-Rem] recently found some examples where the
affine lamination is not locally compact.

As $\mc{R}_c$ is locally compact, we can consider the one point
compactification $\h{\mc{R}}_{c}$ of $\mc{R}_c$; let $*$ be the
point at infinity.  A path $\gamma:[0,1)\rightarrow \mc{R}_c$
\textit{escapes to infinity} in $\mc{R}_c$ if it eventually leaves
every compact set $K\subset \mc{R}_c$. Equivalently, $\gamma$
escapes to infinity if admits an extension
$\h\gamma:[0,1]\rightarrow \h{\mc{R}}_c$ with $\h{\gamma}(1)=*$.
Two paths, $\gamma_1$ and $\gamma_2$, escaping to infinity are
\textit{homotopic at infinity} if for every compact set $K\subset
\mc{R}_c$ there is an $r\geq0$ such that the sub-paths
$\gamma_1|_r:[r,1)\rightarrow \mc{R}_c$ and
$\gamma_2|_r:[r,1)\rightarrow \mc{R}_c$ are homotopic in
$\mc{R}_c\setminus K$.

\begin{definition} Given a leaf $L$ in $\mc{R}_c$, let $E(L)$ denote
the number of paths non-homotopic at infinity.\end{definition}

Now, for a convex cocompact parameter $c$, since the laminated
Julia set $\mc{J}_c$ is compact, we can describe the number
unbounded Fatou components in $\mc{R}_c$ from a topological point
of view:

\begin{lemma}\label{unbounded-homlink}
If $c$ is a convex compact parameter, then for every leaf $L$ in
$\mc{R}_c$, the number $E(L)$ is equal to the number of unbounded
Fatou components in $L$.
\end{lemma}

\begin{proof} This is a direct consequence of the fact, due to
Lyubich and Minsky, that when $c$ is a convex cocompact parameter
the laminated Julia set $\mc{J}_c$ is compact in $\mc{R}_c$, so
every path escaping to infinity must leave eventually the Julia
set $\mc{J}_c$ and escape through an unbounded Fatou component.
\end{proof}

Let $LU_n$ be the set of leaves with exactly $n$ unbounded Fatou
components.

\begin{corollary}
The cardinality of $LU_n$ is a topological invariant.
\end{corollary}

\begin{proof} The number of non-homotopic paths escaping to infinity
is a topological invariant. So, if $h:\mc{R}_c \rightarrow
\mc{R}_{c'}$ is a homeomorphism between regular parts, then
$E(L)=E(h(L))$ for every leaf $L$ in $\mc{R}_c$.
\end{proof}

\begin{corollary}\label{homlink.1} Let $f_c$ be a convex cocompact
quadratic polynomial such that $c$ does not belong to the interior
of the Main Cardioid, then $E(L(\h{\beta}))=1$.
\end{corollary}

\begin{proof} The $\beta$ fixed point is, by definition, a repelling
fixed point where exactly $1$ external ray lands. So, by
Proposition~\ref{unbounded}, the leaf $L(\h{\beta})$ has exactly
one unbounded Fatou component.
\end{proof}

\begin{corollary}\label{discon.cases} If $f_c$ is a convex cocompact
quadratic polynomial such that $c$ does not belong to the Main
Cardioid, then $\h{\infty}$ is the only disconnectivity point of
$\mc{N}_{c}$.
\end{corollary}

\begin{proof} By going through external rays, all leaves in the
regular part have at least one access to $\h{\infty}$. But leaves
with one unbounded Fatou component have access only to
$\h{\infty}$. By Corollary~\ref{homlink.1}, at least the leaf
corresponding to the $\beta$ fixed point has only one unbounded
Fatou component. So, $\h{\infty}$ is the only irregular point that
can be accessed from from every leaf in $\mc{N}_c$.
\end{proof}

Thus, the solenoidal cone at infinity can be characterized as the
only solenoidal cone in $\mc{N}_c$ which connects all of the
leaves of the regular part in the natural extension.

\begin{lemma}\label{periodic.leaves} If $c$ is convex cocompact,
then every periodic leaf has exactly one periodic point.
\end{lemma}

\begin{proof} It is clear that every periodic point in $\mc{R}_{c}$
lies in a periodic leaf. Now, let $L$ be a periodic leaf in
$\mc{R}_c$ of period $n$, so $\hat{f}_c^n(L)=L$. Because, in the
case of convex cocompact parameters, the affine part coincides
with the regular part, for every leaf $L$ in the regular part
there is a uniformization $\psi:L\ra \C$ which conjugates the map
$\h{f}_c^n:L \rightarrow L$ to an affine map $\hat{f}_c^n(z)=az+b$
where $a$ is a complex number. We claim that $a\neq 1$. If on the
contrary $a=1$, as the Julia set $\mc{J}_c$ is compact, there is a
finite covering of small flow boxes, with the property that the
derivative of $\phi$ is bounded away from zero on the Julia set
$\mc{J}_c$.

Let $\h{z}\in\pi^{-1}( A(\infty))$ be any point on the lift of the
basin of infinity. Take $W$ around $z_0=\pi(\h{z})$ as in the
Shrinking Lemma and let $W'$ be the plaque containing $\h{z}$ in
the fiber of $W$. By assumption, $\h{f}_c^{-n m}(W')$ has the same
diameter for all $m$, since translations are isometries. On the
other hand the diameters of $f_c^{-n m}(W)$ are shrinking to $0$
by the Shrinking Lemma. Moreover, for every neighborhood $V$
around the dynamical Julia set $J(f_c)$ we have
$f_c^{-nm}(W)\subset V$ for large enough $n$. This means that the
derivative of $\pi$, under uniformization $\phi$, shrinks to $0$
which is a contradiction. So, the map $\h{f}_c^n$ can not be
conjugated to a translation in $L$. Therefore $a\neq 1$, which
implies the existence of a periodic point in $L$.
\end{proof}

At every periodic point in the dynamic root cycle there are at
least $2$ landing external rays. Also, by definition, any point in
the dynamic root cycle is on the boundary of at least $1$ Fatou
component. Let $L$ be one of the leaves containing a periodic
point in the lift of the dynamic root cycle. Then, by
Proposition~\ref{unbounded}, $E(L)\geq 3$.

If $c$ is the center of a primitive hyperbolic component, then any
leaf containing a periodic point in the dynamic root cycle has $3$
unbounded Fatou components. However, two of the unbounded Fatou
components are associated to the solenoidal cone at infinity,
whereas the other is associated to a Fatou component in the basin
of attraction of the critical cycle. By
Proposition~\ref{unbounded} and Lemma~\ref{unbounded-homlink}, we
obtain the following characterization of regular parts of
primitive superattracting parameters:

\begin{corollary}
A superattracting parameter $c$ is primitive if and only if there
is a leaf $L$ in $\mc{R}_c$ such that $E(L)=3$, and the classes of
paths non-homotopic at infinity in $L$ belong to more than one
solenoidal cones.
\end{corollary}

In order to make the previous corollary more precise, let us
introduce the concept of ends of laminated sets. Let $\mc{L}$ be a
locally compact laminated set, and consider the one point
compactification $\h{\mc{L}}=\mc{L}\cup \{*\}$. Two paths,
$\sigma:[0,1)\rightarrow \mc{L}$ and $\tau:[0,1)\rightarrow
\mc{L}$ escaping to infinity, are said to be \textit{equivalent at
infinity} if for every compact set $K\subset \mc{L}$ there is a
number $r>0$ such that $\sigma([r,1))$ and $\tau([r,1))$ belong to
the same connected component of $\mc{L}\setminus K$. This is an
equivalence relationship in the set of paths escaping to infinity.

\begin{definition} An \textit{end} of a locally compact laminated
set $\mc{L}$, is an equivalence class of the relationship above
described. Let $End(\mc{L})$ denote the set of ends of $\mc{L}$,
then $\mc{L}\cup End(\mc{L})$ is the \textit{end compactification}
of $\mc{L}$.\end{definition}

By definition, each end contains the homotopic class of each of
its elements, so  equivalence  at infinity is a weaker
relationship than the equivalence relationship of being  homotopic
at infinity.

\begin{lemma}\label{irregular.ends} If the Julia set $J(f_c)$ is
locally connected, for every irregular point $\h{I}$ in $\mc{N}_c$
there is an end in $End(\mc{R}_c)$ associated to $\h{I}$.
\end{lemma}

\begin{proof} Since the coordinates of $\h{I}$ belong to the
postcritical set, the coordinates of $\h{I}$ belong either to the
Julia set or to an attracting or superattracting cycle. In any
case, by the local connectivity of $J(f_c)$, there is a point
$z_0$ in the Fatou set $F(f_c)$ and a path $\gamma$ from $z_0$ to
$i_0=\pi(\h{I})$ such that the trajectory of the path $\gamma$
intersects either the Julia set or the postcritical set exactly at
$i_0$. However, the pullbacks $\{\gamma_n\}$ of $\gamma$ are well
defined in $[0,1)$, and altogether define a path
$\h{\gamma}:[0,1)\rightarrow \mc{R}_c$ that escapes to infinity in
the regular part. Hence, we associate the irregular point $\h{I}$
with the end $[\h{\gamma}]$. Let $\h{\gamma}'$ be any other path
defined as $\h{\gamma}$, we want to check that $\h{\gamma}'$ is
equivalent at infinity with $\h{\gamma}$. By definition,
$\h{\gamma}$ and $\h{\gamma}'$ extend to paths from $[0,1]$ to
$\mc{N}_c$ satisfying $\h{\gamma}(1)=\h{\gamma}'(1)=\h{I}$, so the
trajectories of $\h{\gamma}$ and $\h{\gamma}'$ eventually belong
to any neighborhood of $\h{I}$ in $\mc{N}_c$. Since for every
$t\in [0,1)$ the points $\h{\gamma}(t)$ and $\h{\gamma}'(t)$
belong to the regular part, for every compact set $K\subset
\mc{R}_c$ the paths $\gamma$ and $\gamma'$ eventually belong to
the same connected component in $\mc{R}_c\setminus K$.
\end{proof}

\begin{lemma}\label{vert.sol.con}
Let $c$ be a parameter in the Mandelbrot set. There is one and
only one end $E_{\h{\infty}}\in End(\mc{R}_c)$ associated to
$\h{\infty}$ in $\mc{R}_c$. Furthermore, if $c$ is
superattracting. Then every end of $\mc{R}_c$ is associated to a
unique irregular point.
\end{lemma}

\begin{proof} Consider the equipotential $E_r$ in the dynamical
plane. Since $E_r$ is compact and does not contain postcritical
points, the corresponding solenoidal equipotential
$\mc{S}_r=\pi^{-1}(E_r)$ is compact in $\mc{R}_c$. Now, let
$R_\theta$ be any external ray in the dynamical plane, then the
end $[\h{R}]$ of any lift $\h{R}$ of $R_\theta$ in $\mc{R}_c$ is
associated to $\h{\infty}$. If $\h{\gamma}$ is equivalent at
infinity to $\h{R}$ then $\h{\gamma}$ must eventually lie in the
same connected component of $\mc{R}_c\setminus \mc{S}_r$ as
$\h{R}$. This implies that $\pi(\h{\gamma})$ converges to $\infty$
in the dynamical plane, and so $\h{\gamma}$ must converge to
$\h{\infty}$ in $\mc{N}_c$. If $c$ is superattracting, then
instead of the equipotential $E_r$, we can consider an internal
equipotential inside the corresponding Fatou component on the
basin of attraction of the critical cycle. The argument goes
through using the corresponding internal solenoidal equipotential.
\end{proof}

\begin{lemma}\label{sup.irre.end}
For any quadratic polynomial $f_c$, every end of $\mc{R}_c$ is
associated to an irregular point in $\mc{N}_c$.
\end{lemma}

\begin{proof} Let $[\gamma]$ be an end of $\mc{R}_c$, with $\gamma$ a
representative of this end in $\mc{R}_c$. Let $A_n$ denote the
accumulation set of $\gamma_n=\pi_{-n}(\gamma$). Let us check that
$f_c(A_n)=A_{n-1}$. By continuity, $f_c(A_n)\subset A_{n-1}$. Now
let $y\in A_{n-1}$. There is a sequence $t_m \nearrow 1$ such that
$\gamma_{n-1}(t_m)$ converges to $y$. That means that
$\gamma_n(t_m)$ is as close as we want to a point in
$f_c^{-1}(y)$. Since $f_c$ has finite degree, $\gamma_n(t_m)$ must
actually converge to a point in $f_c^{-1}(y)$. So,
$f_c(A_n)=A_{n-1}$ as we claimed, and this implies that we can
construct a backward orbit $\h{y}\in \mc{N}_c$, such that
$\pi_{-n}(\h{y})\in A_n$. Let us check that $\h{y}$ must be
irregular. If, on the contrary, $\h{y}$ is regular, then there is
a $N$ such that $y_n$ is outside the postcritical set of $f_c$ for
$n>N$. Since the postcritical set is closed, there is a
neighborhood $U$ of $y_n$ such that $U$ is outside $P(f_c)$. Let
$K\subset U$ be a compact neighborhood of $y_n$. Since
$\pi_{-n}^{-1}$ is proper, the set $\pi_{-n}^{-1}(K)$ is a compact
neighborhood of $\h{y}$, but $\gamma_n(t_m)$ converges to $y_n$,
so $\gamma(t_m)$ is contained in $K$ for $m$ large. This
contradicts the fact that $\gamma$ is escaping to
infinity.\end{proof}

Lemma~\ref{irregular.ends} does not rule out the possibility that
several irregular points are associated to the same end. It is not
clear whether there is a one-to-one correspondence between the
irregular points and the ends of the regular part. This would
imply that the natural extension of every quadratic parameter $c$
corresponds to the end compactification of $\mc{R}_c$. However, we
have a positive answer for certain parameters.

\begin{proposition}\label{irregular.ends.prop} Let $c$ be a parameter such
that $J(f_c)$ is locally connected, and $c$ is either convex
cocompact or the postcritical set of $f_c$ is a Cantor set. Then
the set of ends corresponds to the set of irregular points, and
the end compactification of $\mc{R}_{c}$ is $\mc{N}_{c}$.
\end{proposition}

\begin{figure}[htbp]
  \begin{center}
  \includegraphics[width=3 in]{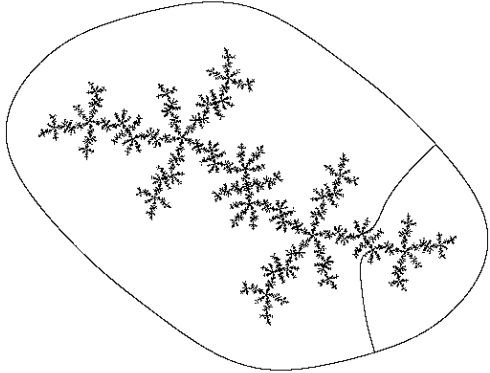}
  \caption{To the proof of Proposition~\ref{irregular.ends.prop}.}\label{cantorproof}
  \end{center}
\end{figure}

\begin{proof} If $f_c$ is convex cocompact then, by Lyubich and
Minsky, the Julia set $\mc{J}_c$ is compact. Hence, the only
irregular points in $\mc{N}_c$ correspond to attracting cycles. By
Proposition~\ref{attract.sup} and Lemma~\ref{vert.sol.con},
attracting cycles correspond to vertices of solenoidal cones, and
vertices of solenoidal cones correspond to ends of the regular
part. Thus, if $c$ is convex cocompact, the end compactification
of $\mc{R}_c$ is homeomorphic to the natural extension $\mc{N}_c$.

Now, assume that the postcritical set is a Cantor set. In this
case, there are no bounded Fatou components in the dynamical
plane. By Lemma~\ref{irregular.ends} and Lemma~\ref{sup.irre.end}
it is enough to prove that different irregular points are
associated to different ends. Let $\h{I}$ and $\h{I}'$ be two
irregular points in $\mc{N}_c$. Without loss of generality we can
assume that $\pi(I)=i_0\neq i_0'=\pi(I')$, also by
Lemma~\ref{irregular.ends} assume that both $i_0$ and $i'_0$
belong to the Julia set $J(f_c)$.

By the local connectivity of $J(f_c)$, there is a path $\gamma$
embedded in $J(f_c)$ connecting $i_0$ with $i'_0$. Let $U$ and
$U'$ be neighborhoods around $i_0$ and $i'_0$ small enough that
$\{t\in[0,1]|\gamma(t)\notin U\cup U' \}$ contains an interval
$(t_1,t_2)$. Since $P(f_c)$ is a Cantor set, there is a $t'\in
(t_1,t_2)$ and an open neighborhood $V$ around $\gamma(t')$ and
not intersecting $P(f_c)$. Since $V$ is open, there are two
external rays $R$ and $R'$ landing at both sides of $\gamma$ in
$V$, say at $z_1$ and $z_2$, and such that the path $\tau$
embedded in $J(f_c)$ from $z_1$ to $z_2$ lies completely in $V$.

Let $T$ be the image of $\tau$ in the dynamical plane. By
construction, the curve whose trajectory is $\sigma=R\cup T \cup
R'$ separates $U$ from $V$. Finally, for any equipotential $E_r$,
the set that consists of the union of $E_r$ and the part of
$\sigma$ inside $E_r$ is a compact set $K$. By construction $K$
does not intersect the postcritical set, so $\pi^{-1}(K)$ is a
compact set in $\mc{R}_c$ such that $I$ and $I'$ lies in different
connected components in $\mc{N}_c$ (see Figure~\ref{cantorproof}).
\end{proof}

\begin{corollary}\label{extension} Let $c$ be any parameter as
in Proposition~\ref{irregular.ends.prop}, then every homeomorphism
$h:\mc{R}_{c}\rightarrow\mc{R}_{c'}$ between regular parts extends
to a homeomorphism of the natural extensions
$\tilde{h}:\mc{N}_c\rightarrow \mc{N}_{c'}$. Moreover,
$\tilde{h}(\h{\infty})=\h{\infty}$
\end{corollary}

\begin{proof} By Proposition~\ref{irregular.ends.prop}, the end compactification of $\mc{R}_c$ is
homeomorphic to the natural extension $\mc{N}_c$, since any
homeomorphism $h:\mc{R}_c\rightarrow \mc{R}_{c'}$ extends to the
end compactification. By Corollary~\ref{discon.cases},
$\h{\infty}$ is the only disconnection point among the irregular
points, and the second part of the corollary follows.
\end{proof}

\begin{corollary}\label{same.period} If
$h:\mc{R}_{c}\rightarrow\mc{R}_{c'}$ is a homeomorphism between
the regular parts of two superattracting parameters $c$ and $c'$,
then the periods of $c$ and $c'$ are equal.
\end{corollary}

\begin{proof} By Corollary~\ref{extension}, $h:\mc{R}_c\rightarrow
\mc{R}_{c'}$ extends to a homeomorphism of the natural extensions
sending irregular points into irregular point. If $p$ is the
period of $c$, there are $p+1$ irregular points in $\mc{N}_c$.
\end{proof}

\begin{lemma}\label{dyn.root.reg}
Let $c_1$ and $c_2$ be two superattracting parameters, and
$h:\mc{R}_{c_1}\rightarrow \mc{R}_{c_2}$ be a homeomorphism. Then
$h$ sends the leaves containing the dynamic root cycle of $c_1$
into the leaves of the dynamic root cycle of $c_2$.
\end{lemma}

\begin{proof} By Proposition~\ref{unbounded} there are at least
three unbounded Fatou components associated to the leaves
containing the dynamic root cycle. On the other hand, these
unbounded Fatou components correspond to at least two different
ends in the regular part. So, one of the ends is associated to at
least two unbounded Fatou components. This end must be
$\h{\infty}$ and no other periodic point in the regular part can
have access to different ends.
\end{proof}

Given a superattracting parameter $c$, let $v_c$ be the valence of
the dynamic root point.

\begin{corollary}
Let $c_1$ and $c_2$ be two superattracting parameters with the
same period. If $v_{c_1}\neq v_{c_2}$, then the corresponding
regular parts $R_{c_1}$ and $R_{c_2}$ are not homeomorphic.
\end{corollary}

\begin{proof} Having different valence, the corresponding leaves
containing the lifts of the dynamic root cycles on each regular
part must have different numbers of unbounded Fatou components,
which implies that the corresponding leaves have different number
of accesses to infinity.
\end{proof}

 \section{Isotopies}
In the previous section we found some invariants that a
homeomorphism between regular spaces must have. In order, to
finish the proof of the Main Theorem, we need to discuss the
isotopy classes of homeomorphisms of regular spaces. In this
section, we will see that among the homeomorphisms isotopic to a
given one, there is always a homeomorphism holding special
properties.

A \textit{self-embedding} of a compact topological space $X$ is a
continuous injective map of $X$ into itself. A homeomorphism
$h:X\rightarrow X$ is called \textit{isotopic to the identity} if
there exist a continuous map $\Phi:[0,1]\times X\rightarrow X$
such that $\Phi(0,x)=x$, $\Phi(1,x)=h(x)$ and the restriction
$\Phi_t(x)=\Phi(t,x)$ is a homeomorphism of $X$ for every $t\in
[0,1]$. In general, two maps $h_1:X\rightarrow Y$ and
$h_2:X\rightarrow Y$ are called \textit{isotopic} if there is a
homeomorphism $\phi:Y \rightarrow Y$ isotopic to the identity,
such that $\phi\circ h_1=h_2$. In this section, we will prove that
every self-embedding of the solenoidal cone $\cone$ is isotopic to
a self-embedding of $\cone$ sending $\sol$ to a solenoid of the
form $\sol\times\{r\}$ in $\cone$.

\subsection{Isotopies of the dyadic solenoid}

The group $H(\sol,\sol)$ of self-homeomorphisms of the solenoid,
endowed with the uniform topology, is a topological group. An
\textit{automorphism} of $\sol$ is an element in $H(\sol,\sol)$
that preserves the group structure of $\sol$. The set of
automorphisms of the solenoid is denoted by $Aut(\sol)$. For
$\tau\in \sol$, the map $\zeta \mapsto \tau\cdot \zeta$ is called
a \textit{left translation} of $\sol$. Abusing notation, we
identify the map with the element $\tau\in \sol$. The set
$Aff(\sol)$ of \textit{affine maps} of the solenoid is a
transformation in $\sol \times Aut(\sol)$, where the first factor
corresponds to the set of left translations in $\sol$.

Let $\ell^2$ be the standard Hilbert space. The following result,
due to James Keesling \cite{Kees}, describes the topological
embedding of $Aut(\sol)$ inside $H(\sol,\sol)$.

\begin{proposition}[Keesling] \label{Keesling}
The group $H(\sol,\sol)$ is homeomorphic to $\ell_2 \times \sol
\times Aut(\sol)$.
\end{proposition}

Actually, Keesling's proof shows that $H(\sol,\sol)/Aff(\sol)$ is
homeomorphic to $\ell_2$. Since $\ell_2$ is a vector space, this
implies that every self-homeomorphism of the solenoid is isotopic
to an affine map. According to Jaroslaw Kwapisz \cite{Kwap}, the
group $Aut(\sol)$ has a simple set of generators:

\begin{proposition} [Kwapisz]\label{Kwapisz} The group $Aut(\sol)$
is the infinite dihedral group generated by $\h{f}_0$ and the
inversion  $s\mapsto \bar{s}$.
\end{proposition}

This proposition was proved in \cite{Kwap} in the more general
setting of $P$-adic solenoids, where $P$ is an arbitrary sequence
of prime numbers. Together these propositions yield the following
corollary:

\begin{corollary}[Kwapisz]\label{solenoid.kwap} Every
homeomorphism of the dyadic solenoid onto itself is isotopic to an
affine map of the form $\tau \circ \h{f}^n \circ r$, where $r$ is
the identity if the map is orientation preserving, or the
inversion $s \mapsto \bar{s}$ otherwise.
\end{corollary}

The following is a known topological property of the dyadic
solenoid, see \cite{Kwap}.

\begin{lemma}\label{indecomposable}
The image of very continuous map $\phi:\sol \rightarrow \sol$, of
the solenoid into itself, is either a point, a closed interval or
onto. \end{lemma}

\begin{proof} The solenoid is a connected, compact metric
Haussdorf space. So is $\phi(\sol)$ by continuity. Now consider a
leaf $L \subset \sol$, then its image $\phi(L)$ is contained in a
leaf $L'$. We claim that, if $\phi(L)$ is unbounded in $L'$ then
$\overline{\phi(L)}=\sol$. If the image of $L$ is a complete leaf
then this is clear because of the density of leaves. Now assume
that $\phi(L)$ is a half line in the solenoid. Recall that the
leaf containing the unit is a one parameter subgroup of the
solenoid. By homogeneity, assume that the unit belongs to
$\phi(L)$ and, after identifying with the real numbers
$\mathbb{R}$, that $\phi(L)$ covers the positive numbers. Now let
$-M$ be some negative number. Since the numbers $2^m$
transversally converge to $0$ as $m$ goes to infinity, the numbers
$2^m-M$ converge to $-M$ as $m$ goes to infinity, so the closure
of $\phi(L)$ contains the whole leaf containing the unit and our
claim follows.

Assume that $\phi(L)$ is bounded, then by connectivity
$\phi(\sol)$ is on the connected component of a bounded set,
therefore it should be an interval or a point. \end{proof}

\subsection{Isotopies of solenoidal cones}

Recall that we can regard the solenoid $\sol$ as the quotient of
$S=I\times F$ by the  map $\sigma$, where $\sigma$ is the
generator of the adding machine action. Here, $I=[0,1]$ and $F$
denotes the fiber over $1$ of the projection $\pi:\sol \rightarrow
\mathbb{S}^1$; $F$ is homeomorphic to the Cantor set $\cantor$.
The cylinder $\sol \times I$ also can be expressed as the quotient
of $S\times I$ by the map $\sigma \times Id$. For $x\in F$, let
$R_x=(0,x)\times I \subset S \times I$, then $R=\cup R_x$ over all
$x\in F$ is the vertical section of $\sol \times I$ which is
homeomorphic to the trivial one dimensional lamination $F \times
I$. The set $R$ can be regarded as the fiber of an external ray in
$\cone$ minus the point $\h{\infty}$. The goal of this section is
to prove:

\begin{figure}[htbp]
\begin{center}
\input{solemb2.pstex_t}
\caption{Embedding of $\sol$ into $\sol\times I$.}
 \label{solemb}
\end{center}
\end{figure}

\begin{proposition}\label{isotopy} Let $\phi:\cone \rightarrow \cone$ be
an orientation preserving self-embedding with $\phi(\sol)\cap \sol
=\emptyset$, then $\cone \setminus \phi(\cone)$ is homeomorphic to
$\sol \times I$.
\end{proposition}

First, we say that two 1-dimensional laminations embedded into a
third lamination of dimension 2 \textit{intersect transversally}
if they intersect leafwise transversally. Transversality is a
smooth notion, so we need an isotopy that regularizes the
embedding $\phi$ on $\sol$. The existence of such an isotopy is
given by a generalization of the corresponding theorem about
surfaces. Namely, if $\gamma:I\rightarrow I\times I$ is a curve in
the unit square, such that $\gamma(I)\cap \partial(I\times
I)=\emptyset$, then there is a map $h:I\times I \rightarrow
I\times I$ isotopic to the identity, $rel$ the endpoints of
$\gamma$, such that $h\circ\gamma$ is piecewise linear, and $h$
leaves the extremes of $\gamma$ fixed. See \cite{Eps}.

Moreover, the corresponding map from the set of embeddings
$\textnormal{Emb}(I,I\times I)$ to the set of self-homeomorphisms
of the unit square isotopic to identity can be chosen to depend
continuously on parameters. That is, if $X$ is a topological
space, and the family of maps $\phi_x:I\ra I\times I$ in
$Top(I,I\times I)$ depends continuously on $x$, then there is a
family of maps $h_x$ in $Top(I\times I,I\times I)$ depending
continuously on $x$ such that $h_x\circ \gamma$ is piecewise
linear. So, making $X=F$, we have the following lemma:

\begin{lemma}\label{pl.boxes}
Let $\phi: S \rightarrow S \times I$ be a laminar embedding, such
that $\phi(S)\cap \partial(S\times I)=\emptyset$. Then, there is a
homeomorphism $h:S\times I\rightarrow S\times I$, isotopic to the
identity, such that $h\circ \phi$ is piecewise linear.
\end{lemma}

This immediately implies:

\begin{lemma}\label{pl.sole}
Given an embedding $\phi:\sol \rightarrow \partial (\sol\times I)$
such that $\phi(\sol)\cap \partial (\sol\times I)=\emptyset$,
there exists a homeomorphism $h:\sol\times I\rightarrow \sol
\times I$, isotopic to the identity, such that $h\circ\phi$ is
piecewise linear.
\end{lemma}

\begin{proof} Let $\{B_1,B_2,...,B_k\}$ be a partition by square flow boxes
in $\sol$ such that each plaque of $B_i$ intersects at most one
vertical segment $R_x$. Apply Lemma~\ref{pl.boxes} to each $B_i$.
\end{proof}

By further local isotopies, we can assure that $\phi(\sol)$
\textit{does not contain vertical segments and that $\phi(\sol)$
intersects $R$ transversally.}

\begin{lemma} \label{inter.fin} For every $x \in F$, the
intersection $\phi(\sol)\cap R_x$ consists of a finite number of
points.
\end{lemma}

\begin{proof} By transversality, every intersection is leafwise isolated.
Since every $R_x$ is compact, there are finitely many
intersections in every $R_x$.
\end{proof}

We can identify $R_x$ with $I$ for each $x\in F$.

\begin{lemma}\label{minimum} Let $m(x)=\min\{r | r \in R_x \cap \phi(\sol)\}$.
Then, the function $m:F\rightarrow I$ depends continuously on $x$.
\end{lemma}

\begin{proof} By compactness and Lemma~\ref{inter.fin}, there is a
covering of $\phi(\sol) \cap R$, by flow boxes
$\{B_1,B_2,...,B_k\}$ in $\sol \times I$  and such that each
plaque of $B_i$ contains a single point of $\phi(\sol) \cap R$ for
every $i$.

The intersection $\phi(\sol)\cap R\cap B_i$ is a transversal for
$B_i$. The points where $m$ attains the minimums are arranged into
these transversals.  By definition, transversals depend
continuously on the fiber. \end{proof}

Fix $x\in F$. By the action of $\sigma$, the segments
$\{R_{\sigma^{n}(x)}\}$ with $n\in \mathbb{Z}$ are precisely the
segments in $R$ that belong to the same leaf in $\sol \times I$.
Let $L$ be the corresponding leaf in $\sol$ that contains
$\phi^{-1}(x)$. By a suitable parametrization, identify $L$ with
$\mathbb{R}$ in an order preserving way.

\begin{lemma} \label{timing} Assume that $\phi$ is orientation preserving and
let $t_x=\phi^{-1}((x,m(x)))$, then $t_x<t_{\sigma(x)}$.
\end{lemma}

\begin{figure}[htbp]
\begin{center}
\input{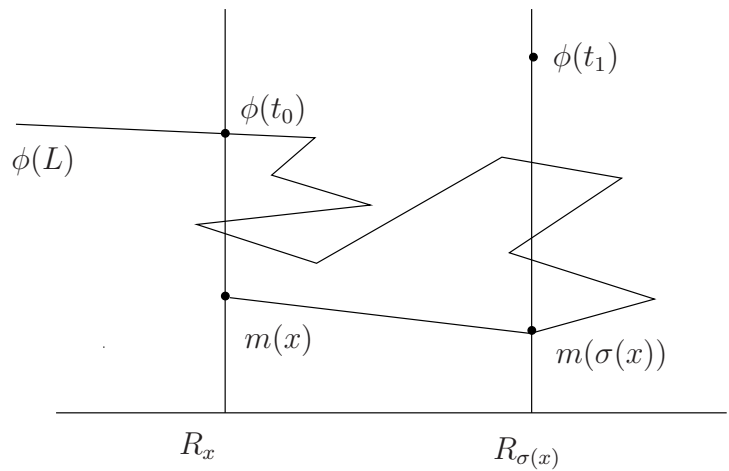}
\caption{To the proof of Lemma~\ref{timing}.}
 \label{timing.fig}
\end{center}
\end{figure}

\begin{proof}  Every leaf $S$ in $\sol\times [0,1]$
is an infinite horizontal strip. By Lemma \ref{indecomposable},
composition of $\phi$ with the vertical projection over the
solenoid is onto. This implies that of every $x\in F$, there is a
first time $t_0$ such that $\phi(t_0)\in R_x$ and a last time
$t_1$ such that $\phi(t_1)\in R_{\sigma(x)}$. Suppose, on the
contrary, that $t_{\sigma(x)}<t_{x}$, since $\phi$ is orientation
preserving, $t_0<t_{\sigma(x)}$, also we have $t_1>t_{\sigma(x)}$.
By definition $\phi(t_0)> m(x)$, and $\phi(t_1)> m(\sigma(x))$ in
$R_x$ and $R_{\sigma(x)}$ respectively. Now, from $m(x)$ there is
no way that the trajectory of $\phi$ gets to $\phi(t_1)$ without
self-intersecting or crossing $R_{\sigma(x)}$ in a lower point
than $m(\sigma(x))$, see Figure~\ref{timing.fig}. Therefore
$t_x<t_{\sigma(x)}$.
\end{proof}

\begin{figure}[htbp]
\begin{center}
\input{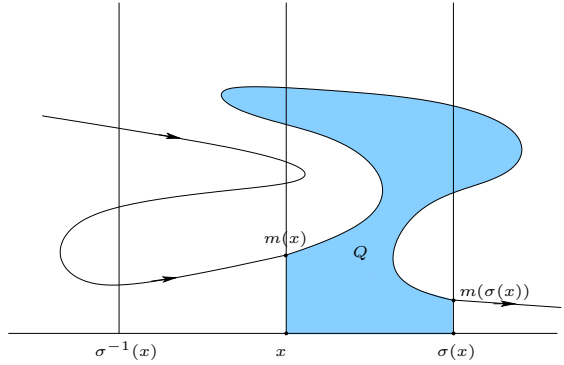}
\caption{The quadrilateral $Q_x$.}
 \label{fig.quadrilateral}
\end{center}
\end{figure}

\begin{proof}[Proof of proposition~\ref{isotopy}.] Let
$(t_x,t_{\sigma(x)})$ be the arc joining $t_x$ with
$t_{\sigma(x)}$ in $\sol$, by definition all of these arcs are
disjoint. Moreover, $\sol=\underset{x\in F}{\cup
}\overline{(t_x,t_{\sigma(x)})}$. By Lemma~\ref{minimum}  the
extremes of these arcs depend continuously on $F$. Now, the arcs
$(m(x),(x,0))$ along $R_x$, also depend continuously on $F$. So,
for every $x\in F$ we have a quadrilateral $Q_x$ in $I\times S$,
with vertices $(m(x),m(\sigma(x)), (\sigma(x),0),(x,0))$, see
Figure~\ref{fig.quadrilateral}. For every $x\in F$, consider
homeomorphisms $\psi_x$ from $Q_x$ to the plaque $I\times I \times
x$ in $I\times S$ such that  $\psi_x=\psi_{\sigma(x)}$ on the side
$(m(\sigma(x)),(\sigma(x),0)$. The homeomorphisms $\psi_x$ paste
together to form the desired homeomorphism from $\cone \setminus
\phi(\cone)$ to $I\times \sol$.\end{proof}

\begin{corollary}\label{frontera1} Let $\phi:\cone \rightarrow \cone$
be an embedding with $\phi(\sol)\cap \sol =\emptyset$, then there
is an isotopy of $\phi$ to a map that sends the boundary onto the
boundary.
\end{corollary}

\begin{proof} The external ray foliation in $\sol \times I$ gives the track of the
desired isotopy.\end{proof}

\begin{lemma}\label{extension.iso} A map
$h:\mc{S}_{R}\rightarrow \mc{S}_{R}$, isotopic to the identity,
extends to a map $\tilde{h}:\cone \rightarrow \cone$ isotopic to
the identity.  This extension can be done in such a way that
$\tilde{h}$ restricted to the complement $\cone \setminus N$ of
some neighborhood $N$ of $S_{R}$ is the identity.
\end{lemma}

\begin{proof} Consider
$N=\cup_{t\in(R-\epsilon,R+\epsilon)}\mc{S}_{t}$, clearly $N$ is a
neighborhood of $\mc{S}_{E_R}$. Let $\Phi: I\times
\mc{S}_{E_T}\rightarrow \mc{S}_{E_R}$ be an isotopy of $h$ to the
identity. Define $b:I\rightarrow I$ by
$b(t)=max\{1,|\frac{\epsilon}{2}(t-R)|\}$. The map
$\tilde{h}:\cone \rightarrow \cone$ given by
$\tilde{h}(t,s)=(t,\Phi(b(t),s))$ satisfies the conditions of the
lemma.\end{proof}

\begin{corollary} Let $\phi:\sol \rightarrow I \times \sol $ be an
embedding. Assume that $\phi$ admits an extension to a map $\sol
\times (-\epsilon,\epsilon) \rightarrow I\times \sol$ for some
positive $\epsilon$. Then, the image of $\phi$ is isotopic to
$\sol \times 0$.
\end{corollary}

\begin{proposition}\label{equipotentials} Assume that $h:\mc{N}_c\rightarrow
\mc{N}_{c'}$ is a homeomorphism such that
$h(\h{\infty})=\h{\infty}$. Then $h$ is isotopic to a
homeomorphism sending a solenoidal equipotential $\mc{S}_R(c)$
onto $\mc{S}_R(c')$.
\end{proposition}

\begin{proof} The solenoidal cone at infinity admits a foliation by
solenoidal equipotentials. These solenoidal equipotentials define
a local basis of neighborhoods homeomorphic to solenoidal cones
around $\h{\infty}$. Locally, $h$ must send this solenoidal local
basis of neighborhoods into a local basis of neighborhoods around
$\h{\infty}$ in $\mc{N}_{c'}$. Thus,  there is a solenoidal cone
at $\mc{S}_R(c)$ in $\mc{N}_c$ which is embedded by $h$ into some
solenoidal cone at $\mc{S}_R(c')$ in $\mc{N}_{c'}$. Now, the
proposition follows from Proposition \ref{isotopy}.
\end{proof}

\subsection{Proof of the Main theorem}

In this section, we will prove the central theorem of this work.
The idea is to recognize the topological imprints that
combinatorics of parameters impose over the regular parts. We will
prove the following:

\begin{theorem}\label{main} Let $h:\mc{R}_{f_c}\rightarrow \mc{R}_{f_{c'}}$ be
an orientation preserving homeomorphism between the regular parts
of two superattracting quadratic polynomials, $f_c$ and $f_{c'}$.
Then, $f_c=f_{c'}$.
\end{theorem}

\begin{proof} By Corollary~\ref{extension}, the homeomorphism $h$
admits an extension $\tilde{h}:\mc{N}_{c}\rightarrow \mc{N}_{c'}$
such that $\tilde{h}(\h{\infty})=\h{\infty}$. By Proposition
\ref{equipotentials}, $\tilde{h}$ is isotopic to a homeomorphism
that sends a solenoidal equipotential $\mc{S}_R$ in $\mc{N}_c$
onto a solenoidal equipotential $\mc{S}_R'$ in $\mc{N}_{c'}$.
Using the lift of B\"ottcher's coordinate at the solenoidal cones
at infinity in $\mc{N}_c$ and $\mc{N}_{c'}$ the map $\tilde{h}$
restricted to $\mc{S}_R$, becomes a homeomorphism of $\sol$ into
itself. By, Corollary~\ref{solenoid.kwap} there exists a map, say
$\psi:\sol \ra \sol$, isotopic to the identity such that $\psi
\circ \tilde {h}$ is an affine transformation of $\sol$ of the
form $\tau\circ \h{f}^n_0$. By Lemma~\ref{extension.iso}, $\psi$
extends to a map $\tilde{\psi}$ isotopic to the identity, defined
on a neighborhood $N$ of $\mc{S}_R$, so that $\tilde{\psi}\circ
\tilde{h}$ coincides with $\tilde{h}$ outside $N$. The map
$\h{f}^{-n}_{c'}\circ\tilde{\psi}\circ\tilde{h}$ is conjugate to
$\tau$ in $\mc{S}_R$. So, by means of the previous normalizations
we can assume that $\tilde{h}$ restricted to $\mc{S}_R$ is already
the translation $\tau$. On the other hand, by
Lemma~\ref{dyn.root.reg}, $\tilde{h}$ must send the leaves
containing the dynamic root cycle of $\h{f}_c$ into the leaves
containing the dynamic root cycle of $\h{f}_{c'}$. The lift of the
dynamic root cycle of $f_c$  may not be mapped to the lift of the
dynamic root cycle of $f_{c'}$. However, on the level of
solenoidal equipotentials, $\tilde{h}$ maps the periodic leaves in
the solenoid $\mc{S}_R$ (under doubling) associated to the ray
portrait of $r_c$, into the corresponding periodic leaves of
$\mc{S}'_R$. Remind that there is a one-to-one correspondence
between periodic leaves and periodic points on $\sol$. Thus,
possibly after another isotopic deformation of $\tilde{h}$, we can
assume that $\tilde{h}|_{\mc{S}_R}$ sends a periodic point in
$\sol$ to a periodic point in $\sol$. But if a translation $\tau$
in $\sol$ sends a periodic point of $\sol$ into a periodic point
in $\sol$, then $\tau$ itself must be periodic. By
Corollary~\ref{same.period}, the periods of $c$ and $c'$ must be
the same.

Therefore, $\tau$ sends every periodic point in $\mc{S}_R\cap
L(\h{r}_c)$ into every periodic point in $\mc{S}'_R\cap
L(\h{r}_{c'})$. That happens in every leaf containing the lift of
the dynamic root cycle. So, by projecting onto $\mathbb{S}^1$ by
$\pi$, the action of $\tau$ becomes a rotation in $\mathbb{S}^1$
that sends the ray portrait of $r_c$ onto the ray portrait of
$r_c'$. By Lemma~\ref{orbit}, the dynamic root cycle must be the
same and then $c=c'$. \end{proof}

We can extend the previous theorem to the following:

\begin{theorem} Let $c$ be a convex cocompact parameter with $Im(c)\neq 0$,
if $h:\mc{R}_c \rightarrow \mc{R}_c'$ is an orientation preserving
homeomorphism between regular spaces then $c=c'$.
\end{theorem}

\begin{proof}
If $c$ is hyperbolic, by Proposition~\ref{attract.sup} the regular
part $\mc{R}_c$ is homeomorphic to $\mc{R}_{c_0}$, where $c_0$ is
the center of the hyperbolic component containing $c$, so we are
in the same situation of Theorem~\ref{main}.

If $c$ is not hyperbolic, then $c\in \partial M$ and since $c$ is
not a real number, $f_c$ has a cycle with valence bigger than 2.

Since the critical point belongs to the Julia set and there are no
irregular points in $\fiber(J(f))$, the point $\h{\infty}$ is the
only irregular point in $\mc{N}_c$. Moreover, the laminated Julia
set $\mc{J}_c$ is compact, so it follows that the end
compactification of $\mc{R}_c$ is homeomorphic to $\mc{N}_c$.
Because $h$ is a homeomorphism, the end compactification of
$\mc{R}_{c'}$ also coincides with $\mc{N}_{c'}$, and $\mc{N}_{c'}$
has only one irregular point, thus $c'$ is also a convex cocompact
parameter. Moreover, $h$ extends to a homeomorphism sending
$\h{\infty}_c$ to $\h{\infty}_{c'}$.

We are now in the situation of the hypothesis of
Proposition~\ref{equipotentials}, and we can find a homeomorphism
$\tilde{h}$ in the isotopy class of $h$, such that $\tilde{h}$
sends a solenoidal equipotential in $\mc{R}_c$ onto another in
$\mc{R}_{c'}$. When we restrict to these solenoidal
equipotentials, the map acts as a translation $\tau$.

By the remark above, there is a cycle with valence bigger than 2,
so the lift of this cycle belongs to a cycle of leaves with at
least 3 unbounded Fatou components. The image of any of these
leaves is also a leaf $L$ with at least 3 unbounded Fatou
components; by Proposition~\ref{unbounded}, $L$ is periodic under
$f_{c'}$. Using a similar argument as in the proof of
Theorem~\ref{main}, we can check that $\tilde{h}$ must preserve
the ray portraits associated to the cycle of $L$. Moreover, the
translation $\tau$ is equal to the identity. Hence, preserves the
orbit portrait of every cycle of $f_c$. The same holds for
$\tilde{h}^{-1}$, then $c$ and $c'$ are combinatorially
equivalent. But if $c$ is convex cocompact, by Carleson, Jones and
Yoccoz \cite{CJY}, $c$ is rigid which implies $c=c'$.
\end{proof}

When $c$ is real and convex cocompact, all periodic leaves of
$\mc{R}_c$ have either 1 or 2 unbounded components. However, there
are many non-periodic leaves with 2 unbounded components. So, the
argument of unbounded components does not guarantee that $h$ must
send periodic leaves into periodic leaves, and there is no reason
for the ray portraits to be preserved. In this case, we would need
another argument to justify that the induced translation on
solenoidal equipotentials is the identity.

 \nocite{*}

\bibliographystyle{amsplain}
\bibliography{class}

\end{document}